# UN THÉORÈME LIMITE POUR LES COVARIANCES DES SPINS DANS LE MODÈLE DE SHERRINGTON–KIRKPATRICK AVEC CHAMP EXTERNE

By Albert Hanen

*Université Paris X*

On étudie la covariance (pour la mesure de Gibbs) des spins en deux sites dans le cas d'un modèle de Sherrington–Kirkpatrick avec champ externe; lorsque le nombre de sites du modèle tend vers l'infini, une évaluation asymptotique des moments d'ordre $p$ de cette covariance permet d'obtenir un théorème limite faible avec une loi limite en général non gaussienne.

We study the covariance (for Gibbs measure) of spins at two sites in the case of a Sherrington–Kirkpatrick model with an external field. When the number of sites of the model grows to infinity, an asymptotic evaluation of the $p$ moments of that covariance allows us to obtain a weak limit theorem, with a generally non-Gaussian limit law.

**1. Introduction.** Nous nous placerons dans le cadre d'un modèle particulier de verres de spins, celui de Sherrington–Kirkpatrick (SK) avec champ externe (voir [9]), tel qu'il est étudié par Talagrand au Chapitre 2 de son livre [10]. Donnons quelques définitions et notations: $N$ étant un entier, on appelle configuration (ou $N$-configuration) de spins toute suite

$$\sigma = (\sigma_1, \ldots, \sigma_i, \ldots, \sigma_N) \in \Sigma_N \stackrel{\mathrm{def}}{=} \{-1, +1\}^N.$$

La variable $\sigma_i$ est le spin au site $i$ de la configuration $\sigma$.

Soient $\beta$ et $h$ deux réels positifs, représentant respectivement l'inverse d'une température et un champ magnétique externe constant. Soit $g_{i,j}$, $i<j$, une famille de variables i.i.d. gaussiennes standards représentant le désordre du modèle. On appelle hamiltonien du modèle l'application qui à toute $N$-configuration $\sigma$ associe le nombre

$$(1) \qquad -H_N(\sigma, \beta, h) = \frac{\beta}{\sqrt{N}} \sum_{1 \leq < i < j \leq N} g_{i,j} \sigma_i \sigma_j + h\left(\sum_{i \leq N} \sigma_i\right).$$

Considérons la distribution des masses $\exp(-H_N(\sigma,\beta,h))$ sur l'espace $\Sigma_N$. La masse totale portée par $\Sigma_N$ est appelée la fonction de partition du modèle et se note $Z_N$. La mesure de Gibbs est la probabilité sur $\Sigma_N$, notée $G_N$, telle que

$$G_N(\sigma) = \exp(-H_N(\sigma,\beta,h))/Z_N.$$

On introduit également le recouvrement de deux configurations $\sigma^1$ et $\sigma^2$, qui est le nombre

(2) $$R(\sigma^1,\sigma^2); \quad \text{noté aussi} \quad R_{1,2} = \frac{\sum_{i=1}^{i=N} \sigma_i^1 \sigma_i^2}{N}.$$

Lorsque le nombre de sites $N \to +\infty$, les propriétés asymptotiques de la fonction de partition $Z_N$, étudiées via celles de la variable $Y_N \stackrel{\text{def}}{=} \log(Z_N)/N$, ainsi que celles des recouvrements, jouent un rôle très important dans la compréhension du modèle à haute température, c'est-à-dire pour $\beta$ suffisamment petit, de même que pour des modèles voisins (voir [9]). Ces propriétés sont maintenant connues pour l'essentiel.

Dans le cas du modèle SK avec champ externe, en se plaçant à haute température, Talagrand montre en particulier, dans le Chapitre 2 de son livre [10], la convergence de $E(Y_N)$ si $N \to +\infty$ et calcule sa limite (c'est la "solution réplique symétrique"); il montre aussi la constance asymptotique, en un certain sens, des recouvrements et établit divers théorèmes centraux limites concernant leurs fluctuations; par ailleurs, un théorème central limite (TCL) concernant les fluctuations de $Y_N$ est démontré dans l'article [6].

Ces résultats ont été précédés d'études de même type dans le cas $h=0$. En particulier, dans l'article [5], les auteurs introduisent des outils de calcul stochastique qui leurs permettent de retrouver un théorème limite pour les fluctuations de $Y_N$ démontré dans l'article [1] et montrent un TCL pour les lois des recouvrements construits sur $n$ répliques, $\Sigma_N^{\otimes n}$ étant muni de la probabilité $G_N^{\otimes n}$.

Les fluctuations de $Y_N$ et des recouvrements ont été également étudiées dans le cas des $p$-spins; il s'agit de modèles dans lesquels les sous ensembles $J$ de $p$ sites engendrent des interactions $g_J$ i.i.d. gaussiennes standards. On peut alors construire l'analogue de l'hamiltonien (1), et bâtir une mesure de Gibbs correspondante.

Dans le cas $h=0$, Talagrand, au Chapitre 6 de son livre [10], construit une constante $\gamma_p$ telle que $E(Y_N)$ converge pour $\beta < \gamma_p$ si $N \to +\infty$, définissant ainsi une zone à haute température.

Toujours dans le cadre des modèles de $p$-spins sans champ externe, dans le travail [4], les auteurs, s'appuyant sur les outils de calcul stochastique développés dans l'article [5], montrent la normalité asymptotique des fluctuations de $Y_N$ sur une échelle polynomiale dépendant de $p$, et ce pour $\beta < \tilde{\gamma}_p$, où $\tilde{\gamma}_p < \gamma_p$.

Kurkova, dans l'article [8], étend ce résultat à toute la zone à haute température définie par Talagrand ($\beta < \gamma_p$), et établit dans les mêmes conditions un TCL pour les recouvrements, avec la même échelle de fluctuations que pour $p = 2$.

Enfin, dans le travail [2], les auteurs étudient le modèle à $p$ spins avec champ externe et montrent, en utilisant ou étendant les méthodes de Talagrand, l'analogue des résultats obtenus à haute température pour $p = 2$ dans son livre [10].

Nous nous intéresserons quant à nous à l'étude du comportement asymptotique de la covariance, pour la mesure de Gibbs $G_N$, des spins en deux sites $i$ et $j$, quantité notée $\gamma_{i,j}$. Cette quantité est liée aux recouvrements et l'on ne trouve dans la littérature que peu de résultats la concernant; le résultat asymtotique le plus précis est donné par [10] (Corollaire 2.6.2), qui permet d'obtenir à haute température la limite du moment d'ordre 2 de $N^{1/2}\gamma_{i,j}$ si $N \to +\infty$.

Nous étudierons le comportement asymptotique des moments d'ordre $p$ de $N^{1/2}\gamma_{i,j}$. Soit $z$ une variable gaussienne standard et soit

(3) $\quad Y = \beta z\sqrt{q_2} + h,\qquad$ où $q_2$ vérifie la relation $E(th^2(Y)) = q_2$.

On sait (voir [10], Proposition 2.4.8) que $q_2$ (notée $q$ dans [10]) existe et est unique si $h > 0$.

Soit $U = 1 - th^2(Y)$. Soient $z_1$ et $z_2$ deux variables gaussiennes centrées réduites, indépendantes entre elles et indépendantes de $z$, engendrant des variables de type $U$, notées $U_1$ et $U_2$, i.i.d. de loi celle de $U$. Nous montrerons le théorème suivant:

THÉORÈME 1.1. *Pour $\beta$ suffisamment petit, les moments d'ordre $p$ de $N^{1/2}\gamma_{i,j}$ convergent, si $N \to +\infty$, vers ceux de la variable $\frac{\beta}{\sqrt{1-\beta^2 E(U^2)}} z U_1 U_2$.*

REMARQUE 1.2. En particulier, pour le moment d'ordre 2, on retrouve la limite $\sigma^2 \stackrel{\text{def}}{=} \frac{\beta^2}{1-\beta^2 E(U^2)} E^2(U^2)$, donnée par Talagrand. D'autre part, la variable limite est centrée, mais non gaussienne, car ses moments d'ordre $2q$ ne sont pas de la forme $\sigma^{2q} E(z^{2q})$. Toutefois, lorsque $h = 0$ et $\beta < 1$, on a $q_2 = 0, U = 1$. La loi limite est alors gaussienne, résultat déjà annoncé par Talagrand ([10], Research Problem 2.3.10).

*Notation $O(k)$* (*voir* [10], 2.99)  Nous utiliserons constamment la notation définie çi après suivante:

DEFINITION 1.3. Nous dirons qu'une expression numérique, notée $F(N,\beta,\theta)$, définie pour $\beta \leq \beta_0$, vérifie la relation

$$F(N,\beta,\theta) = O(k),$$

s'il existe une constante $K$, indépendante de $\beta$ et de $N$, mais pouvant dépendre du paramètre $\theta$, telle que

$$|F(N,\beta,\theta)| \leq \frac{K}{N^{k/2}}.$$

Introduisons d'autres notations utilisées par Talagrand dans le cadre de notre modèle: Si $f$ est une fonction de $k$ $N$-configurations, appelées aussi répliques, $\langle f \rangle$ est l'intégrale de f par rapport à la probabilité produit $G_N^{\otimes k}$. La mesure de Gibbs étant aléatoire, $\langle f \rangle$ est aussi aléatoire, nous poserons

$$\nu(f) = E(\langle f \rangle).$$

Remarquons que $\langle f \rangle$, et donc aussi $\nu(f)$, sont invariantes par permutation des répliques en jeu. Remarquons également que les mesures $\nu$ (mais pas les mesures produits $G_N^{\otimes k}$) sont invariantes par permutation des $N$ sites.

Nous montrerons que l'étude des moments d'ordre $p$ de la covariance se ramène à l'étude de $\nu(f)$, pour un certain $f$ dépendant de $2p$ répliques. Pour calculer ou évaluer de telles expressions, nous emploierons la méthode dite du "smart path," développée par Talagrand tout au long de son livre [10]. Celle ci consiste pour l'essentiel à construire une famille continue d'hamiltoniens $-H_{N,t}(\sigma,\beta,h)$, $t \in (0,1)$ (l'expression exacte de ces hamiltoniens est donnée en Annexe), engendrant de la même manière chacun une mesure de Gibbs, une intégrale et son espérance que nous noterons naturellement respectivement: $G_{N,t}, \langle f \rangle_t$ et $\nu_t(f) = E(\langle f \rangle_t)$, elles aussi invariantes par permutation des répliques en jeu. Le cas $t = 1$ correspond à l'hamiltonien (1), le cas $t = 0$ à une situation plus simple, soit içi pour notre modèle, celle où, pour la mesure de Gibbs $G_{N,0}$, les spins aux $N-1$ premiers sites sont indépendants du spin au site $N$. $\nu_t(f)$ est dérivable (de tous ordres) en $t$ et cette propriété peut permettre le calcul de $\nu(f)$ par un développement de Taylor de $\nu_t(f)$, pour $t = 1$, en $t = 0$.

La méthode du "smart path" est au coeur des études d'autres modèles de verres de spins réalisées par Talagrand dans son livre [10] pour résoudre le même type de problématiques que celles envisagées dans le modèle SK, à savoir établir une "solution réplique symétrique" et des TCL pour les recouvrements. Il s'agit entre autres, outre le cas SK avec $h = 0$, de modèles de réseaux de neurones comme celui du perceptron, avec des interactions gaussiennes ou de Bernouilli, des configurations de spins éventuellement à valeurs dans une sphère de $\mathbf{R}^N$, ou du modèle de Hopfield avec champ externe.

Revenons à notre modèle pour observer qu'en général, le développement de Taylor de $\nu_t(f)$ à l'ordre 1 suffit à Talagrand pour obtenir ses résultats. Signalons toutefois l'article [3], utilisant également cette méthode, où des développements à un ordre $> 1$ sont donnés dans le cas $h = 0$ pour étudier les moments d'ordre $n$ des recouvrements pour la mesure $\nu$.

Dans notre cas, pour l'étude des moments d'ordre $p$ de la covariance, un développement de Taylor à l'ordre $p$ sera nécessaire. Nous montrerons que pour l'étude asymptotique, on peut se ramener au cas $p$ pair et à l'évaluation du terme principal du développement $\nu_0^{(p/2)}(f)/(p/2)!$, pour laquelle nous utiliserons notamment un TCL sur les recouvrements.

REMARQUE 1.4. Dans la suite, on considèrera souvent plusieurs répliques, indexées supérieurement et notées généralement $\sigma^r$, dont les spins au site $i$ seront notés $\sigma_i^r$ (ou $\varepsilon^r$ si $i = N$).

La preuve détaillée du résultat principal sera donnée dans la suite de ce travail, qui est organisée de la manière suivante: La Section 2 donne une formule permettant, $f$ étant une fonction numérique de plusieurs répliques, de donner la forme des dérivées de $\nu_t(f)$ d'ordre $l > 1$ et d'obtenir leur ordre de grandeur asymptotique.

La section suivante, intitulée Outils, étudie des expressions du type $\nu_t(f^- \tilde{\varepsilon}^{\otimes p})$, où $f^-$ ne dépend pas des spins au site $N$ et

$$(4) \qquad \tilde{\varepsilon}^r = \varepsilon^{2r-1} - \varepsilon^{2r}, \qquad \tilde{\varepsilon}^{\otimes p} = \prod_{r=1}^{r=p} \tilde{\varepsilon}^r.$$

On montre pour l'essentiel, grâce à une proposition de nature combinatoire, que les dérivées d'ordre $l$ en $t = 0$ de $\nu_t(f^- \tilde{\varepsilon}^{\otimes p})$ sont nulles si $2l < p$, et l'on en donne l'expression exacte lorsque $2l = p$.

La section qui suit est divisée en deux sous sections; dans la première, il est d'abord montré que l'étude des moments d'ordre $p$ de la covariance $\gamma_{i,j}$ se ramène à celle de $\nu_1(f^- \tilde{\varepsilon}^{\otimes p})$, où

$$(5) \qquad \text{si } \tilde{\sigma}_i^r = \sigma_i^{2r-1} - \sigma_i^{2r}, \qquad f^- = \tilde{\sigma}_1^{\otimes p} = \prod_{i=1}^{i=p} \tilde{\sigma}_1^r.$$

On montre ensuite que, dans ce cas particulier, l'expression de la dérivée d'ordre $l = p/2$ en $t = 0$ de $\nu_t(f^- \tilde{\varepsilon}^{\otimes p})$ obtenue dans la section précédente se simplifie.

Dans la seconde sous section, il est montré que les dérivées d'ordre $l$ en $t = 0$ de $\nu_t(f^- \tilde{\varepsilon}^{\otimes p})$ sont des $O(p+1)$ lorsque $2l > p$. Ceci permet, par un développement de Taylor à l'ordre $p$, de se ramener au cas où $p$ est pair ($p = 2q$) et à l'évaluation de la dérivée d'ordre $q$ de $\nu_t(f^- \tilde{\varepsilon}^{\otimes 2q})$ en $t = 0$. Cette étude fait l'objet de la section qui suit. En étudiant l'expression simplifiée de cette dérivée obtenue plus haut, et en utilisant notamment un résultat obtenu dans la section "Outils" ainsi qu' un TCL pour les recouvrements, on parvient à une évaluation asymptotique précise de cette dérivée.

Nous pouvons alors obtenir, dans la dernière partie, la loi limite, au sens de la convergence des moments, des variables $N^{1/2}\gamma_{i,j}$, annoncée dans le Théorème 1.1.

Ce travail est complété par deux annexes. L'Annexe A expose de manière plus détaillée les méthodes de la cavité et du "smart path," puis développe quelques notions étudiées dans le livre [10] qui interviennent de manière essentielle dans ce travail, notamment les inégalités exponentielles et les TCL sur les recouvrements.

L'Annexe B est consacrée à la démonstration de résultats techniques énoncés dans les deux parties consacrées respectivement au développement de Taylor à l'ordre $p$ de $\nu_t(f^-\tilde{\varepsilon}^{\otimes p})$ et à l'étude des moments d'ordre pair de la covariance.

REMARQUE 1.5. Dans toute la suite de ce travail, nous nous placerons dans l'hypothèse $\beta \leq \beta_0$, où $\beta_0$ est suffisamment petit pour que les inégalités exponentielles et les TCL sur les recouvrements évoqués dans l'Annexe A soient valables.

Les résultats démontrés dans ce travail ont été annoncés dans la note [7].

## 2. Dérivées d'ordre $\geq 1$ de $\nu_t(f)$.

REMARQUE 2.1. A l'exception de l'intervalle $[0,1]$ de **R**, les intervalles considérés dans la suite seront des intervalles d'entiers.

Soit $f$ une fonction de $n$ répliques $\sigma^1, \sigma^2, \ldots, \sigma^n$. Talagrand ([10], Proposition 2.4.5) donne une formule permettant de calculer la dérivée de $\nu_t(f)$. Introduisons les notations suivantes, s'ajoutant à celles des recouvrements définis dans l'Introduction [voir la relation (2)]:

NOTATIONS 2.2. Posons
$$R_{1,2}^- = R^-(\sigma^1, \sigma^2) = \frac{\sum_{i=1}^{i=N-1} \sigma_i^1 \sigma_i^2}{N},$$
$$R_{1,2}' = \frac{N}{N-1} R_{1,2}^-.$$

(Cette dernière quantité est le recouvrement des configurations concernées restreintes aux $N-1$ premiers sites.)

- On pose $\dot{R}_J = R_{l_1,l_2} - q_2 = R(\sigma^{l_1}, \sigma^{l_2}) - q_2$ où $J = \{l_1, l_2\}$ $(1 \leq l_1 < l_2)$, $q_2$ vérifie la relation (3).
- On pose $\dot{R}_J^- = R^-(\sigma^{l_1}, \sigma^{l_2}) - q_2$.
- On pose $\varepsilon^J = \varepsilon^{l_1} \varepsilon^{l_2} (\varepsilon^r = \sigma_N^r)$.
- On note $D^1(n)$ l'ensemble des sous ensembles à deux éléments ou doublets $\{l_1, l_2\}, 1 \leq l_1 < l_2 \leq n$.

Posons également

$$c_J(n) = \begin{cases} 1, & \text{si } J \in D^1(n), \\ -n, & \text{si } J \text{ a pour extrémité droite } n+1, \\ n(n+1)/2, & \text{si } J = \{n+1, n+2\}, \\ 0, & \text{si } J = \{l_1, n+2\}, \text{ où } l_1 \leq n. \end{cases}$$

On a la formule suivante, qui est une réécriture de l'équation (2.170) du livre [10]:

(6) $$\nu'_t(f) = \beta^2 \sum_{J \in D^1(n+2)} c_J(n) \nu_t(f \varepsilon^J \dot{R}^-_J).$$

On voit que l'argument de chacun des $\nu_t$ du membre de droite dépend de $n$ (respectivement $n+1$, $n+2$) répliques si $c_J(n) = 1$ [respectivement $c_J(n) = -n, c_J(n) = n(n+1)/2$ ou $0$]. Cette formule se prête à itérations et l'on peut donner, au prix de notations parfois lourdes, l'expression exacte des dérivées d'ordre $> 1$ de $\nu_t(f)$.

Nous démontrerons le théorème suivant, qui suffira pour la suite:

THÉORÈME 2.3. *Soit $f$ une fonction numérique de $n$ répliques. Soit $l \geq 1$ un entier, $\hat{J} = (J_1, \ldots, J_l)$ une suite quelconque de doublets de l'intervalle d'entiers $(1, 2, \ldots, n+2l)$. La dérivée d'ordre $l \geq 1$ de $\nu_t(f)$ est donnée par la formule*

(7) $$\nu_t^{(l)}(f) = \beta^{2l} \sum_{\hat{J}} c_{\hat{J}}(n) \nu_t\left(f \prod_{i=1}^{i=l} (\varepsilon^{J_i} \dot{R}^-_{J_i})\right).$$

*De plus, les coefficients $c_{\hat{J}}(n)$ peuvent être éventuellement nuls, mais sont égaux à $1$ pour toutes les suites $\hat{J}$ dont chaque terme $J_i \in D^1(n)$.*

On en déduira l'évaluation asymptotique suivante en utilisant la notation $O(k)$ définie dans l'Introduction (Définition 1.3) et qui sera constamment utilisée dans la suite:

COROLLAIRE 2.4. *Pour tout $t \in [0,1]$,*

$$|\nu_t^{(l)}(f)| = \nu^{1/2}(f^2) O(l).$$

DÉMONSTRATION DU THÉORÈME 2.3. Elle se fait par récurrence sur l'ordre de dérivation $l$. La formule (7), ainsi que les conditions sur les coefficients, sont vraies pour $l = 1$, d'après l'équation (6). Supposons les vraies à l'ordre $l$; pour chaque suite $\hat{J}$ apparaissant dans l'écriture de $\nu_t^{(l)}(f)$, il suffira de dériver par la formule (6) le terme $\nu_t(f \prod_{i=1}^{i=l}(\varepsilon^{J_i} \dot{R}^-_{J_i})) \stackrel{\text{def}}{=} \nu_t(g_1)$.

La fonction $g_1$ dépend de $n_{\hat{j}}$ répliques avec $n \le n_{\hat{j}} \le n + 2l$. La dérivée fera apparaitre une combinaison linéaire de termes $\beta^2 \nu_t(g \varepsilon^{J_{l+1}} \dot{R}^-_{J_{l+1}})$, avec un certain coefficient valant 1 si $J_{l+1} \in D^1(n) \subset D^1(n_{\hat{j}})$. En combinant le résultat obtenu avec l'écriture de $\nu_t^{(l)}(f)$, on obtient la formule (7) au rang $l+1$, ainsi que les propriétés des coefficients. $\square$

DÉMONSTRATION DU COROLLAIRE 2.4. Soit
$$g = \prod_{i=1}^{i=l} (\varepsilon^{J_i} \dot{R}^-_{J_i}) \quad \text{alors} \quad |g| = \prod_{i=1}^{i=l} |\dot{R}^-_{J_i}|.$$

On a, en utilisant l'équation (72),
$$\left| \nu_t \left( f \prod_{i=1}^{i=l} (\varepsilon^{J_i} \dot{R}^-_{J_i}) \right) \right| \le \nu_t(|f||g|) \le K \nu(|f||g|).$$

L'inégalité de Schwarz montre que
$$\langle |f||g| \rangle \le \langle f^2 \rangle^{1/2} \langle g^2 \rangle^{1/2},$$
$$E(\langle |f||g| \rangle) \le E(\langle f^2 \rangle^{1/2} \langle g^2 \rangle^{1/2}) \le E^{1/2}(\langle f^2 \rangle) E^{1/2}(\langle g^2 \rangle),$$

d'où
$$\nu(|f||g|) \le \nu^{1/2}(f^2) \nu^{1/2}(g^2).$$

Les inégalités de Hölder montrent alors que
$$\langle g^2 \rangle = \left\langle \prod_{i=1}^{i=l} (\dot{R}^-_{J_i})^2 \right\rangle \le \prod_{i=1}^{i=l} \langle (\dot{R}^-_{J_i})^{2l} \rangle^{1/l}.$$

D'où
$$\nu(g^2) = E(\langle g^2 \rangle) \le E\left( \prod_{i=1}^{i=l} (\langle (\dot{R}^-_{J_i})^{2l} \rangle^{1/l}) \right) \le \prod_{i=1}^{i=l} E^{1/l}(\langle (\dot{R}^-_{J_i})^{2l} \rangle)$$

et donc
$$\nu^{1/2}(g^2) \le \prod_{i=1}^{i=l} \nu^{1/(2l)}((\dot{R}^-_{J_i})^{2l}).$$

On sait, d'après l'inégalité exponentielle (77), que
(8) $\quad E^{1/(2l)}(\langle (\dot{R}^-_J)^{2l} \rangle) = \nu^{1/(2l)}[(\dot{R}^-_J)^{2l}] = O(1).$

Donc $\nu^{1/2}(g^2) = O(l)$. Le Théorème 2.3 montre que $\nu_t^{(l)}(f)$ est combinaison linéaire de termes du type $\nu_t(fg)$ qui sont chacun des $O(l)\nu^{1/2}(f^2)$, et donc est aussi un $O(l)\nu^{1/2}(f^2)$. $\square$

**3. Outils.** Soit $f^-$ une fonction numérique définie sur $\sum_{N-1}^{\otimes n'}$, $p$ un entier, $n = \sup(n', 2p)$. Supposons que l'on se donne une suite de $n$ répliques $\sigma^1, \ldots, \sigma^n$, le spin au site $N$ de la réplique $\sigma^r$ étant noté $\varepsilon^r$, d'après la Remarque 1.4.

Nous nous intéresserons au calcul d'expressions du type $\nu(f^- \tilde{\varepsilon}^{\otimes p})$, où $\tilde{\varepsilon}^{\otimes p}$ est défini par la relation (4), à l'aide d'un développement de Taylor de $\nu_t(f^- \tilde{\varepsilon}^{\otimes p})$, pour $t = 1$, en $t = 0$. Cela implique de caractériser les dérivées d'ordre $l$ de $\nu_t(f^- \tilde{\varepsilon}^{\otimes p})$ en $t = 0$.

NOTATIONS 3.1. Nous noterons $\hat{\pi}_q$ une partition de l'intervalle $(1, 2, \ldots, 2q)$ en $q$ doublets. On peut ordonner entre eux ces doublets, obtenant ainsi une partition ordonnée de $(1, 2, \ldots, 2q)$, notée $\pi_q^*$. Il y a $q!$ partitions ordonnées distinctes associées à une même partition $\hat{\pi}_q$, qui sera dite non-ordonnée. Le $i^{\text{eme}}$ doublet d'une partition ordonnée $\pi_q^*$ sera noté $\{r_i, s_i\}$, où $r_i < s_i$.

REMARQUE 3.2. Nous identifierons parfois une partition non ordonnée à la partition ordonnée dite naturelle obtenue en supposant la suite $r_i$ croissante.

Posons
$$\tilde{\sigma}^r = \sigma^{2r-1} - \sigma^{2r}, \quad \text{c'est-à-dire} \quad \tilde{\sigma}_i^r = \sigma_i^{2r-1} - \sigma_i^{2r} \qquad \forall i \leq N.$$

Soient
$$R(\tilde{\sigma}^r, \tilde{\sigma}^s) = \frac{\sum_{j=1}^{j=N} \tilde{\sigma}_j^r \tilde{\sigma}_j^s}{N},$$
$$R^-(\tilde{\sigma}^r, \tilde{\sigma}^s) = \frac{\sum_{j=1}^{j=N-1} \tilde{\sigma}_j^r \tilde{\sigma}_j^s}{N}.$$

On notera
$$R_{\pi_q^*} = \prod_{i=1}^{i=q} R(\tilde{\sigma}^{r_i}, \tilde{\sigma}^{s_i}),$$
$$R_{\pi_q^*}^- = \prod_{i=1}^{i=q} R^-(\tilde{\sigma}^{r_i}, \tilde{\sigma}^{s_i}).$$

Ces quantités sont indépendantes de l'ordre des doublets et pourront être notées également respectivement: $R_{\hat{\pi}_q}$, $R_{\hat{\pi}_q}^-$.

Nous montrerons le théorème suivant:

THÉORÈME 3.3. *Soit $f^-$ définie sur $\sum_{N-1}^{n'}$. Soit $A_q' = E(\frac{1}{ch^{4q}(Y)}) = E(1 - th^2(Y))^{2q}$. On a la relation*

(9) $\quad \nu_0^{(l)}(f^- \tilde{\varepsilon}^{\otimes p}) = \begin{cases} 0, & \text{si } 2l < p, \\ q! \beta^{2q} A_q' \sum_{\hat{\pi}_q} \nu_0(f^- R_{\hat{\pi}_q}^-), & \text{si } 2l = p = 2q. \end{cases}$

*On a alors*

$$\nu(f^-\tilde{\varepsilon}^{\otimes 2q}) = \frac{\nu_0^{(q)}(f^-\tilde{\varepsilon}^{\otimes 2q})}{q!} + O(q+1)\nu^{1/2}((f^-)^2). \tag{10}$$

Nous en déduirons les deux corollaires suivants, qui seront aussi utilisés dans les sections suivantes:

COROLLAIRE 3.4. *Soit $C$ un ensemble d'entiers disjoint de $(1, 2, \ldots, 2p)$. Soit $\eta_C = \{\varepsilon^r : r \in C\}$ et $h_1$ une fonction des variables de $\eta_C$, c'est-à-dire une fonction des spins au site $N$ de configurations indexées par $C$. On a*

$$\nu_0^{(q)}(f^-\tilde{\varepsilon}^{\otimes 2q}h_1(\eta_C)) = q!\beta^{2q}E\left(\frac{\langle h_1(\eta_C)\rangle_0}{ch^{4q}(Y)}\right)\sum_{\hat{\pi}_q}\nu_0(f^-R^-_{\hat{\pi}_q}), \tag{11}$$

$$\nu(f^-\tilde{\varepsilon}^{\otimes 2q}h_1(\eta_C)) = \beta^{2q}E\left(\frac{\langle h_1(\eta_C)\rangle_0}{ch^{4q}(Y)}\right)\sum_{\hat{\pi}_q}\nu_0(f^-R^-_{\hat{\pi}_q})$$
$$+ O(q+1)\nu^{1/2}((f^-)^2). \tag{12}$$

Introduisons les notations suivantes, qui seront utilisées dans la suite:

- Soit $B$ un sous ensemble de l'intervalle $(1, 2, \ldots, n)$, de cardinal $|B|$.
- Posons $\varepsilon^B = \prod_{r\in B}\varepsilon^r$, $r(B) = \sum_{x\in B}(x+1)$.

On a le corollaire suivant:

COROLLAIRE 3.5. *Soit $B'$ un sous ensemble d'entiers tel que $|B'| \leq 2k < p$. Soit $q$ tel que $p = 2q$ ou $2q - 1$. On a la relation*

$$\nu_0^{(l)}(f^-\tilde{\varepsilon}^{\otimes p}\varepsilon^{B'}) = \begin{cases} 0, & \text{si } 2l + 2k < p, \\ O(q-k)\nu^{1/2}((f^-)^2), & \text{si } 2l + 2k \geq p. \end{cases} \tag{13}$$

La démonstration du théorème, puis de ses corollaires, passe par celle d'une proposition de nature combinatoire que nous énoncerons après avoir introduit les notations suivantes:

Soit:

- Le symbole $a_r$ dénote l'un des deux nombres $2r-1$ ou $2r$, $a'_r$ l'autre.
- Le symbole $\mathcal{C}$ dénote la classe des sous ensembles $B$ de $(1, 2, \ldots, 2p)$, appelés aussi ensembles $p$-canoniques, s'écrivant sous la forme
$$B = \{a_1, a_2, \ldots, a_p\}.$$

On voit que $|\mathcal{C}| = 2^p$. On sait que l'on a, d'après la relation (4),

$$\tilde{\varepsilon}^{\otimes p} = \prod_{r=1}^{r=p}(\varepsilon^{2r-1} - \varepsilon^{2r}).$$

En développant ce produit comme somme de termes du type $\prod_{r=1}^{r=p}(-1)^{(a_r+1)} \times \varepsilon^{a_r}$, on obtient

$$\tilde{\varepsilon}^{\otimes p} = \sum_{B \in \mathcal{C}} (-1)^{r(B)} \varepsilon^B. \tag{14}$$

PROPOSITION 3.6. *Soient $\widehat{B}$ un sous ensemble de l'intervalle $(1, 2, \ldots, 2p)$ et $C$ un sous ensemble de $(1, 2, \ldots, n)$ disjoint de $(1, 2, \ldots, 2p)$. On a*

$$\langle \tilde{\varepsilon}^{\otimes p} \varepsilon^{\widehat{B}} \varepsilon^C \rangle_0 = \begin{cases} 0, & \text{si } \widehat{B} \notin \mathcal{C}, \\ (-1)^{r(\widehat{B})} \dfrac{(th(Y))^{|C|}}{(ch^2(Y))^p}, & \text{si } \widehat{B} \in \mathcal{C}. \end{cases} \tag{15}$$

On en déduira le corollaire suivant:

COROLLAIRE 3.7. *Plaçons nous dans les conditions du Corollaire 3.4. On a de même*

$$\langle \tilde{\varepsilon}^{\otimes p} \varepsilon^{\widehat{B}} h_1(\eta_C) \rangle_0 = \begin{cases} 0, & \text{si } \widehat{B} \notin \mathcal{C}, \\ (-1)^{r(\widehat{B})} \dfrac{\langle h_1(\eta_C) \rangle_0}{(ch^2(Y))^p}, & \text{si } \widehat{B} \in \mathcal{C}. \end{cases} \tag{16}$$

Nous admettrons provisoirement la Proposition 3.6 et le Corollaire 3.7, qui seront démontrés à la fin de la section.

REMARQUE 3.8. Nous utiliserons dans toute la suite les propriétés de la mesure de Gibbs $G_{N,0}$ décrites dans l'Annexe A, à savoir l'indépendance du vecteur des spins aux $N-1$ premiers sites, qui ont pour loi les mesures de Gibbs $G_{N-1}^-$ (avec les notations afférentes A.1), et du spin au site $N \sigma_N \stackrel{\text{def}}{=} \varepsilon$, qui vérifie la relation (69): $\langle \varepsilon \rangle_0 = th(Y)$. Ces propriétés d'indépendance s'étendent à plusieurs répliques et sont régies par les équations (70) et (71) de l'Annexe A.

Avant de passer aux démonstrations, nous pouvons faire la remarque technique suivante:

REMARQUE 3.9. Si $\triangle$ désigne la différence symétrique ensembliste, on a, si $B_1$ et $B_2$ sont deux ensembles d'entiers indexant des répliques,

$$\varepsilon^{B_1} \varepsilon^{B_2} = \varepsilon^{B_1 \triangle B_2}$$

[car si $r \in B_1 \cap B_2, \varepsilon^r$ intervient deux fois dans $\varepsilon^{B_1} \varepsilon^{B_2}$ et $(\varepsilon^r)^2 = 1$]. La différence symétrique étant associative, on a une relation analogue avec plusieurs ensembles $B_i$.

DÉMONSTRATION DU THÉORÈME 3.3. Soit $g = f^{-}\tilde{\varepsilon}^{\otimes p}$; $g$ dépend de $n = \max(2p, n')$ répliques. On a, en utilisant les relations (70) et (71),

$$\nu_0(g) = \nu_0(f^{-})\nu_0(\tilde{\varepsilon}^{\otimes p})$$
$$= 0 \text{ car } \nu_0(\tilde{\varepsilon}^{\otimes p}) = 0.$$
(17)

En effet, la relation (69) induit $\langle \tilde{\varepsilon}^r \rangle_0 = 0 \; \forall r$ et l'indépendance des répliques pour les mesures produits associées aux intégrales $\langle \cdot \rangle_0$ entraine $\langle \tilde{\varepsilon}^{\otimes p} \rangle_0 = 0$, d'où le résultat (17) en prenant l'espérance.

En développant $\nu_0^{(l)}(f^{-}\tilde{\varepsilon}^{\otimes p})$ selon le Théorème 2.3 précédent, on obtient une combinaison linéaire de termes du type

$$\nu_0(f^{-}\tilde{\varepsilon}^{\otimes p}\varepsilon^{J_1}\cdots\varepsilon^{J_l}\dot{R}_{J_1}^{-}\cdots\dot{R}_{J_l}^{-}).$$

Les $J_i$ sont des doublets du type $\{l_i, l_i'\}$, $l_i < l_i'$. Plus précisément, suivant l'équation (7) écrite pour $t = 0$, en posant $\hat{J} = (J_1, \ldots, J_l)$, on a

$$\nu_0^{(l)}(f^{-}\tilde{\varepsilon}^{\otimes p}) = \beta^{2l} \sum_{\hat{J}} c_{\hat{J}}(n) \nu_0\left(f^{-}\tilde{\varepsilon}^{\otimes p}\prod_{i=1}^{i=l}(\varepsilon^{J_i}\dot{R}_{J_i}^{-})\right).$$
(18)

On a de même, en utilisant les relations (70) et (71),

$$\nu_0\left(f^{-}\tilde{\varepsilon}^{\otimes p}\prod_{i=1}^{i=l}(\varepsilon^{J_i}\dot{R}_{J_i}^{-})\right) = \nu_0(f^{-}\dot{R}_{J_1}^{-}\cdots\dot{R}_{J_l}^{-})\nu_0(\tilde{\varepsilon}^{\otimes p}\varepsilon^{J_1}\cdots\varepsilon^{J_l}).$$
(19)

On écrit, en utilisant la Remarque 3.9, par associativité de la différence symétrique: $\varepsilon^{J_1}\cdots\varepsilon^{J_l} = \varepsilon^{J_1 \triangle \cdots \triangle J_l}$.

On pose

$$(J_1 \triangle \cdots \triangle J_l) \cap (1, 2, \ldots, 2p) = \widehat{B}, (J_1 \triangle \cdots \triangle J_l) \cap (1, 2, \ldots, 2p)^c = C.$$

Alors $\nu_0(\tilde{\varepsilon}^{\otimes p}\varepsilon^{J_1}\cdots\varepsilon^{J_l}) = \nu_0(\tilde{\varepsilon}^{\otimes p}\varepsilon^{\widehat{B}}\varepsilon^C)$. On a aussi

$$|J_1 \triangle \cdots \triangle J_l| \leq |J_1 \cup \cdots \cup J_l| \leq 2l,$$

les égalités n'ayant lieu que si les $J_i$ sont disjoints.

Il est donc clair que $|\widehat{B}| \leq 2l$.

(a) Si $2l < p$, $\widehat{B}$ a moins de $p$ éléments et donc $\widehat{B} \notin \mathcal{C}$. Par conséquent, en vertu de la Proposition 3.6, on a

$$\forall \hat{J} \quad \nu_0(\tilde{\varepsilon}^{\otimes p}\varepsilon^{J_1}\cdots\varepsilon^{J_l}) = \nu_0(\tilde{\varepsilon}^{\otimes p}\varepsilon^{\widehat{B}}\varepsilon^C) = E(\langle\tilde{\varepsilon}^{\otimes p}\varepsilon^{\widehat{B}}\varepsilon^C\rangle_0) = 0.$$

D'où

$$\nu_0^{(l)}(f^{-}\tilde{\varepsilon}^{\otimes p}) = 0.$$

(b) Si $p = 2q$, le même raisonnement montre que l'on a l'équivalence

$$\nu_0(\tilde{\varepsilon}^{\otimes 2q} \varepsilon^{J_1} \cdots \varepsilon^{J_q}) \neq 0$$

(20) $\quad\Longleftrightarrow\quad C = \varnothing$ et les doublets $J_i$ sont deux à deux disjoints, d'union $\widehat{B} \in \mathcal{C}$.

On a dans ce cas, en vertu de la Proposition 3.6,

(21) $\quad \nu_0(\tilde{\varepsilon}^{\otimes 2q} \varepsilon^{J_1} \cdots \varepsilon^{J_q}) = (-1)^{r(\widehat{B})} E\left(\frac{1}{ch^{4q}(Y)}\right) = (-1)^{r(\widehat{B})} A'_q.$

Si les $J_i$ vérifient la relation (20), il existe une partition ordonnée de $(1, 2, \ldots, 2q)$ de type $\pi_q^*$ en $q$ doublets notés respectivement $\{r_i, s_i\}$ telle que

$$\text{lorsque } \widehat{B} = \{a_1, \ldots, a_r, \ldots, a_{2q}\}, \qquad J_i = \{a_{r_i}, a_{s_i}\}.$$

On a de plus

(22) $$r(\widehat{B}) = \sum_{i=1}^{i=q} r(J_i).$$

NOTATION.  Nous noterons

$$\Delta_i = \{\{2r_i - 1, 2s_i - 1\}, \{2r_i - 1, 2s_i\}, \{2r_i, 2s_i - 1\}, \{2r_i, 2s_i\}\}.$$

Il est clair que l'ensemble à quatre éléments $\Delta_i$ est le domaine de variation de $J_i$ lorsque $J_i$ est de type $\{a_{r_i}, a_{s_i}\}$.

Remarquons que réciproquement, si l'on se donne une partition quelconque de type $\pi_q^*$, toute suite de $q$ doublets $J_i$ telle que $\forall i \leq q, J_i \in \Delta_i$ vérifie la relation (20).

De plus, comme $\widehat{B} \subset (1, 2, \ldots, 2p) \subset (1, 2, \ldots, n)$, si l'on revient à l'équation (18) avec $l = q$, pour chacune des suites $\hat{J}$ induisant un terme non nul au membre de droite de cette équation, on a

$$\forall i \leq q \qquad J_i \subset (1, 2, \ldots, n).$$

Le coefficient $c_{\hat{j}}(n)$ du terme correspondant vaut donc 1, d'après le Théorème 2.3. On a $\Delta_i \subset D^1(n)$.

On a donc, en utilisant les équations (18), (19), (21) et ( 22),

(23) $\quad \nu_0^{(q)}(f^- \tilde{\varepsilon}^{\otimes 2q}) = \beta^{2q} A'_q \sum_{\pi_q^*} \sum_{J_i \in \Delta_i, \forall i \leq q} (-1)^{r(J_1) + \cdots + r(J_q)} \nu_0(f^- \dot{R}^-_{J_1} \cdots \dot{R}^-_{J_q}).$

Si l'on pose $g_1(J) \stackrel{\text{def}}{=} (-1)^{r(J)}\dot{R}_J^-$, l'équation (23) s'écrit aussi

$$\nu_0^{(q)}(f^-\tilde{\varepsilon}^{\otimes 2q}) = \beta^{2q}A'_q \sum_{\pi_q^*} \sum_{J_i \in \Delta_i, \forall i \leq q} \nu_0\left(f^- \prod_{i=1}^{i=q}(g_1(J_i))\right) \quad (24)$$

$$= \beta^{2q}A'_q \sum_{\pi_q^*} \nu_0\left(f^- \prod_{i=1}^{i=q}\left(\sum_{J_i \in \Delta_i} g_1(J_i)\right)\right). \quad (25)$$

On a par ailleurs

$$R^-(\tilde{\sigma}^r, \tilde{\sigma}^s) = \frac{1}{N}\sum_{j=1}^{j=N-1} \tilde{\sigma}_j^r \tilde{\sigma}_j^s = \frac{1}{N}\sum_{j=1}^{j=N-1}(\sigma_j^{2r-1} - \sigma_j^{2r})(\sigma_j^{2s-1} - \sigma_j^{2s}).$$

En développant ces produits, on obtient finalement

$$R^-(\tilde{\sigma}^r, \tilde{\sigma}^s) = \sum_{a_r, a_s} (-1)^{(a_r+a_s)} \dot{R}^-_{\{a_r, a_s\}}. \quad (26)$$

D'où les relations,

$$\forall 1 \leq i \leq q \qquad R^-(\tilde{\sigma}^{r_i}, \tilde{\sigma}^{s_i}) = \sum_{J_i \in \Delta_i} g_1(J_i). \quad (27)$$

On a par conséquent, en reportant ce résultat dans l'équation (25), la relation

$$\nu_0^{(q)}(f^-\tilde{\varepsilon}^{\otimes 2q}) = \beta^{2q}A'_q \sum_{\pi_q^*} \nu_0\left(f^- \prod_{i=1}^{i=q} R^-(\tilde{\sigma}^{r_i}, \tilde{\sigma}^{s_i})\right)$$

$$= \beta^{2q}A'_q \sum_{\pi_q^*} \nu_0(f^- R_{\pi_q^*}).$$

On a bien le résultat (9) annoncé.

On sait aussi que le développement de Taylor de $\nu_t(f^-\tilde{\varepsilon}^{\otimes 2q})$ en $t=0$ à l'ordre $q$ fait apparaitre un reste se comportant comme la dérivée d'ordre $q+1$ en $t' \leq t$, qui est, d'après le Corollaire 2.4, un $O(q+1)\nu^{1/2}((f^-)^2)$.

En appliquant ce résultat au cas $t=1$, on trouve bien la dernière partie du théorème. □

REMARQUE 3.10. Talagrand avait obtenu dans la Proposition 2.6.3 de son livre [10] les équations (9) et (10) dans le cas $p=2$, c'est-à-dire $q=1$, où tout se simplifie.

REMARQUE 3.11. On voit, grâce au corollaire 2.4, que si $2l > p$, les dérivées $\nu_t^{(l)}(f^-\tilde{\varepsilon}^{\otimes p})$, $\forall t(0 \leq t \leq 1)$ sont en général des $O(l)\nu^{1/2}((f^-)^2)$.

DÉMONSTRATION DU COROLLAIRE 3.4. La démonstration est analogue à celle du théorème qui précède, en utilisant le Corollaire 3.7. □

DÉMONSTRATION DU COROLLAIRE 3.5. L'expression $\nu_0^{(l)}(\tilde{\varepsilon}^{\otimes p}\varepsilon^{B'}f^-)$ est une combinaison linéaire de termes du type

$$\nu_0(\tilde{\varepsilon}^{\otimes p}\varepsilon^{B'}\varepsilon^{J_1}\varepsilon^{J_2}\cdots\varepsilon^{J_l})\nu_0(f^-\dot{R}_{J_1}^-\cdots\dot{R}_{J_l}^-).$$

Il est clair que le premier facteur est nul si $2l < p - 2k$, car si

$$\varepsilon^{B'}\varepsilon^{J_1}\varepsilon^{J_2}\cdots\varepsilon^{J_l} = \varepsilon^{\widehat{B}_1},$$

$\widehat{B}_1$ ne peut avoir qu'au maximum $|B'| + 2l \leq 2k + 2l$ éléments (lorsque les ensembles considérés définissant $\widehat{B}_1$ sont disjoints) et par conséquent ne peut être de la forme (indiquée dans la Proposition 3.6) $\{a_1,\ldots,a_p\}\cup C$.

Il en résulte bien la première partie de l'équation (13).

On sait, d'après le Corollaire 2.4, que $\nu_0^{(l)}(\tilde{\varepsilon}^{\otimes p}\varepsilon^{B'}f^-) = O(l)\nu^{1/2}((f^-)^2)$ et que $O(l) = O(q-k)$ si $l \geq q - k$. Il en résulte bien la seconde partie de l'équation (13). □

DÉMONSTRATION DE LA PROPOSITION 3.6. Nous commencerons par démontrer le lemme suivant:

LEMME 3.12. *Plaçons nous dans les conditions de la Proposition 3.6. Soit $\triangle$ la différence symétrique ensembliste. Pour $B \in \mathcal{C}$, soit*

(28) $$F(B) = (-1)^{r(B)}th(Y)^{|B\triangle\widehat{B}|}.$$

*On a*

(29) $$\langle\tilde{\varepsilon}^{\otimes p}\varepsilon^{\widehat{B}}\varepsilon^C\rangle_0 = \left(\sum_{B\in\mathcal{C}} F(B)\right)(th(Y))^{|C|}.$$

DÉMONSTRATION. Rappelons que $\langle\varepsilon^B\rangle_0 = (th(Y))^{|B|}$ ([10], relation (2.167)). Rappelons (voir Remarque 3.9) que l'on a:

$$\varepsilon^{B_1}\varepsilon^{B_2} = \varepsilon^{B_1\triangle B_2}.$$

La formule (14) et les propriétés d'indépendance des répliques pour la mesure produit liée à $\langle\cdot\rangle_0$ impliquent que

$$\langle\tilde{\varepsilon}^{\otimes p}\varepsilon^{\widehat{B}}\varepsilon^C\rangle_0 = \sum_{B\in\mathcal{C}}(-1)^{r(B)}\langle\varepsilon^B\varepsilon^{\widehat{B}}\varepsilon^C\rangle_0 = \sum_{B\in\mathcal{C}}(-1)^{r(B)}\langle\varepsilon^{B\triangle\widehat{B}}\varepsilon^C\rangle_0$$

$$= \sum_{B\in\mathcal{C}}(-1)^{r(B)}th(Y)^{|B\triangle\widehat{B}|}th(Y)^{|C|}.$$

[On a bien sur $B \triangle \widehat{B} \subset \{1, \ldots, 2p\}$ donc $(B \triangle \widehat{B}) \cap C = \varnothing$.] On obtient bien la relation (29) annoncée. $\square$

Nous allons montrer maintenant le résultat de nature combinatoire suivant:

LEMME 3.13. *Soient $B = \{a_1, \ldots, a_p\} \in \mathcal{C}$, $B_l$ l'ensemble obtenu en substituant dans $B$ $a_l'$ à $a_l$, les propriétés suivantes sont équivalentes:*

(a) $F(B) = -F(B_l)$,
(b) $1_{\widehat{B}}(a_l) = 1_{\widehat{B}}(a_l')$.

DÉMONSTRATION. Soit $B_0$ le complémentaire de $\{a_l\}$ dans $B$. On peut écrire $B = B_0 \cup \{a_l\}$ et $B_l = B_0 \cup \{a_l'\}$. On a $B_l \in \mathcal{C}$, $(B_l)_l = B$ $\forall B \in \mathcal{C}$ et donc

$$\sum_{B \in \mathcal{C}} F(B) = \sum_{B \in \mathcal{C}} F(B_l). \tag{30}$$

Comparons maintenant $F(B)$ et $F(B_l)$. On a

$$r(B) = r(B_0) + a_l + 1, \qquad r(B_l) = r(B_0) + a_l' + 1.$$

Or $a_l$ et $a_l'$ sont de parité opposée, donc également $r(B)$ et $r(B_l)$. On a donc

$$(-1)^{r(B)} = -(-1)^{r(B_l)}. \tag{31}$$

Comparons $|B \triangle \widehat{B}|$ et $|B_l \triangle \widehat{B}|$. Puisque

$$|B \triangle \widehat{B}| = |B| + |\widehat{B}| - 2|B \cap \widehat{B}|$$

et que $|B| = |B_l| = p$, on a

$$|B \triangle \widehat{B}| = |B_l \triangle \widehat{B}| \iff |B \cap \widehat{B}| = |B_l \cap \widehat{B}|.$$

On a

$$|B \cap \widehat{B}| = |B_0 \cap \widehat{B}| + |\{a_l\} \cap \widehat{B}| = |B_0 \cap \widehat{B}| + 1_{\widehat{B}}(a_l)$$

et de même

$$|B_l \cap \widehat{B}| = |B_0 \cap \widehat{B}| + 1_{\widehat{B}}(a_l').$$

Donc $|B \triangle \widehat{B}| = |B_l \triangle \widehat{B}| \iff 1_{\widehat{B}}(a_l) = 1_{\widehat{B}}(a_l')$. L'équation (28) donne l'équivalence annoncée. $\square$

Venons en maintenant à la démonstration de la Proposition 3.6 proprement dite:

(a) Soit il existe un doublet $\{2l_0-1, 2l_0\} \subset \widehat{B}$, alors $1_{\widehat{B}}(a_{l_0}) = 1_{\widehat{B}}(a'_{l_0}) = 1$ $\forall B \in \mathcal{C}$ et d'après le Lemme 3.13, $F(B) = -F(B_{l_0})$, l'équation (30) montre que $\langle \tilde{\varepsilon}^{\otimes p} \varepsilon^{\widehat{B}} \varepsilon^C \rangle_0 = 0$, dans ce cas.

(b) Soit un tel doublet n'existe pas. Alors $\widehat{B}$ a au plus un élément dans chaque doublet et donc $|\widehat{B}| \leq p$.

1. Soit $|\widehat{B}| < p$, alors il existe un doublet noté $\{2l_0-1, 2l_0\}$ disjoint de $\widehat{B}$. Dans ce cas, $1_{\widehat{B}}(a_{l_0}) = 1_{\widehat{B}}(a'_{l_0}) = 0$ $\forall B \in \mathcal{C}$ et l'on a aussi

$$F(B) = -F(B_{l_0}), \qquad \langle \tilde{\varepsilon}^{\otimes p} \varepsilon^{\widehat{B}} \varepsilon^C \rangle_0 = 0.$$

2. Soit $|\widehat{B}| = p$, alors $\widehat{B}$ a exactement un élément dans chaque doublet, et donc $\widehat{B} \in \mathcal{C}$. Nous écrirons

$$\widehat{B} = \{\hat{a}_1, \ldots, \hat{a}_l, \ldots, \hat{a}_p\}.$$

Dans ce cas, pour tout $B \in \mathcal{C}$, si $B = \{a_1, \ldots, a_l, \ldots, a_p\}$ $\forall l \leq p$, un seul des deux nombres $a_l$ ou $a'_l \in \widehat{B}$, donc l'argument précédent ne s'applique pas.

On va calculer $F(B)$ $\forall B \in \mathcal{C}$. Puisque

$$|B \triangle \widehat{B}| = |B| + |\widehat{B}| - 2|B \cap \widehat{B}|$$

et que $|B| = |\widehat{B}| = p$, on a, si $|B \cap \widehat{B}| = k$,

$$|B \triangle \widehat{B}| = 2(p-k)$$

l'exposant de $th(Y)$ dans $F(B)$ est alors $2(p-k)$ pour de tels $B$; reste à les caractériser et à évaluer $(-1)^{r(B)}$.

On peut écrire, pour de tels $B$,

(32)
$$B = (B \cap \widehat{B}) \cup D, \qquad \widehat{B} = (B \cap \widehat{B}) \cup \widehat{D},$$
$$\text{où } D = B \cap \widehat{B}^c, \widehat{D} = \widehat{B} \cap B^c.$$

On a $|\widehat{D}| = |D| = p-k$.

Remarquons que si $a_l \in D, a_l \notin \widehat{B}$, donc $a_l \neq \hat{a}_l$, et par conséquent on a

$a_l = \hat{a}'_l, \qquad a'_l = \hat{a}_l, \qquad$ donc $\hat{a}_l \in \widehat{D}$ et $a_l = \hat{a}'_l \in \widehat{D}'$, c'est-à-dire $D \subset \widehat{D}'$,

où $\widehat{D}'$ désignant l'ensemble obtenu en remplaçant chaque élément $\hat{a}_l$ de $\widehat{D}$ par $\hat{a}'_l$.

Comme $D$ et $\widehat{D}$ ont le même cardinal, $D = \widehat{D}'$; la donnée de $\widehat{B}$ et de $\widehat{B} \cap B$ détermine donc $B$, par les égalités (32).

Il y a donc $C_p^k$ ensembles $B \in \mathcal{C}$ distincts tels que $|B \cap \widehat{B}| = k$. On a également

$$r(B) = r(B \cap \widehat{B}) + r(D), \qquad r(\widehat{B}) = r(B \cap \widehat{B}) + r(\widehat{D}), \text{ d'où}$$

$$(33) \qquad r(B) = r(\widehat{B}) + r(\widehat{D}') - r(\widehat{D}),$$

$$r(\widehat{D}') - r(\widehat{D}) = \sum_{x \in \widehat{D}} ((x'+1) - (x+1)) = \sum_{x \in \widehat{D}} (x' - x).$$

Remarquons que si $x = \hat{a}_r$ et $x \in \{2r-1, 2r\}$,

$$x' - x = \begin{cases} 1, & \text{si } x \text{ est impair,} \\ -1, & \text{si } x \text{ est pair.} \end{cases}$$

Posons $\widehat{D}^+ = \{x : x \in \widehat{D}, x \text{ impair}\}$ et de même $\widehat{D}^- = \{x : x \in \widehat{D}, x \text{ pair}\}$. On a, d'après les relations (33) et (32),

$$r(\widehat{D}') - r(\widehat{D}) = |\widehat{D}^+| - |\widehat{D}^-|, \qquad p - k = |\widehat{D}| = |\widehat{D}^+| + |\widehat{D}^-|.$$

D'où $r(\widehat{D}') - r(\widehat{D}) = p - k - 2|\widehat{D}^-|$ a la parité de $p - k$ et donc les relations (28) et (33) entrainent:

$$\forall B \in \mathcal{C} \text{ tel que } |B \cap \widehat{B}| = k, \qquad F(B) = (-1)^{r(\widehat{B})} (-1)^{p-k} (th(Y))^{2(p-k)},$$

$$\sum_{B \in \mathcal{C}} F(B) = (-1)^{r(\widehat{B})} \sum_{k=0}^{k=p} (-1)^{p-k} C_p^k (th(Y))^{2(p-k)},$$

$$= (-1)^{r(\widehat{B})} (1 - th^2(Y))^p,$$

$$= (-1)^{r(\widehat{B})} \frac{1}{(ch^2(Y))^p}.$$

On trouve bien le résultat annoncé. □

DÉMONSTRATION DU COROLLAIRE 3.7. Il suffit d'utiliser l'indépendance pour les mesures produits associées à $\langle \cdot \rangle_0$ des répliques d'indice (en exposant) $r$ appartenant à des ensembles disjoints. □

**4. Moments d'ordre $p$ des covariances pour la mesure de Gibbs des spins en deux sites et formule de Taylor.** Certains résultats de cette section seront démontrés dans l'Annexe B.

4.1. *Définitions et premiers résultats.* Nous montrerons dans cette sous section que l'étude des moments d'ordre $p$ des covariances des spins en deux sites se ramène à celle d'une expression de la forme $\nu(f^- \tilde{\varepsilon}^{\otimes p})$, pour un $f^-$ particulier. Le Théorème 3.3 s'applique donc et nous obtenons de plus une simplification de la relation (9).

Énonçons plus précisément nos résultats: Soit $\gamma_{i,j} = \langle \sigma_i \sigma_j \rangle - \langle \sigma_i \rangle \langle \sigma_j \rangle$ la covariance des spins aux sites $i$ et $j$. Soit $\tilde{\sigma} = \sigma^1 - \sigma^2$ la différence de deux répliques.

DÉFINITION 4.1. On appelle covariance symétrisée des spins aux sites $i$ et $j$ la quantité
$$\tilde{\gamma}_{i,j} = \langle \tilde{\sigma}_i \tilde{\sigma}_j \rangle.$$

On voit facilement que

(34) $$\gamma_{i,j} = \tfrac{1}{2}\tilde{\gamma}_{i,j}.$$

L'étude des moments d'ordre $p$ de $\gamma_{i,j}$ se ramène donc à celle des moments d'ordre $p$ de $\tilde{\gamma}_{i,j}$, qui sont indépendants du choix de $i$ et $j$. En prenant $i=1, j=N$, on obtient, en utilisant les notations de la section précédente, le théorème suivant:

THÉORÈME 4.2. *Si* $f^- = \tilde{\sigma}_1^{\otimes p} = \tilde{\sigma}_1^1 \cdots \tilde{\sigma}_1^p$ *alors* (i)

(35) $$E(\tilde{\gamma}_{i,j}^p) = E(\langle \tilde{\sigma}_1 \tilde{\sigma}_N \rangle^p) = \nu(f^- \tilde{\varepsilon}^{\otimes p}).$$

(ii) *Le Théorème* 3.3 *vaut pour cet* $f^-$.

*De plus, lorsque* $p=2q$, *la seconde partie de la relation* (9) *se simplifie et devient*:

(36) $$\nu_0^{(q)}(\tilde{\sigma}_1^{\otimes 2q} \tilde{\sigma}_N^{\otimes 2q}) = \beta^{2q} A'_q \frac{(2q)!}{2^q} E[(\langle \tilde{\sigma}_1^1 \tilde{\sigma}_1^2 R^-(\tilde{\sigma}^1, \tilde{\sigma}^2) \rangle_-)^q].$$

REMARQUE 4.3. Les égalités (35) et le Théorème 3.3 [précisé par la relation (36)] suggèrent d'évaluer le moment d'ordre $p$ de $\tilde{\gamma}_{i,j}$ par un développement de Taylor à un ordre adéquat de $\nu_t(f^- \tilde{\varepsilon}^{\otimes p})$ pour $t=1$ en $t=0$. C'est ce à quoi nous allons nous attacher dans la sous section suivante.

DÉMONSTRATION DU THÉORÈME 4.2. On sait que
$$\langle \tilde{\sigma}_1 \tilde{\sigma}_N \rangle^p = \langle \tilde{\sigma}_1^1 \tilde{\sigma}_N^1 \tilde{\sigma}_1^2 \tilde{\sigma}_N^2 \cdots \tilde{\sigma}_1^p \tilde{\sigma}_N^p \rangle$$
$$= \langle \tilde{\sigma}_1^{\otimes p} \tilde{\sigma}_N^{\otimes p} \rangle = \langle \tilde{\sigma}_1^{\otimes p} \tilde{\varepsilon}^{\otimes p} \rangle.$$

D'où, en prenant l'espérance, la relation (35).

Par ailleurs, pour toute partition $\pi_q^*$ de l'intervalle $(1,2,\ldots,2q)$ en $q$ doublets $\{r_i, s_i\}$ (où $i=1,\ldots,q$), on a $\tilde{\sigma}_1^{\otimes p} = \prod_{i=1}^{i=q}(\tilde{\sigma}_1^{r_i} \tilde{\sigma}_1^{s_i})$. Il en résulte, en utilisant aussi les relations (70) et (71) que

$$\left\langle f^- \prod_{i=1}^{i=q} R^-(\tilde{\sigma}^{r_i}, \tilde{\sigma}^{s_i}) \right\rangle_0 = \left\langle f^- \prod_{i=1}^{i=q} R^-(\tilde{\sigma}^{r_i}, \tilde{\sigma}^{s_i}) \right\rangle_-$$
$$= \left\langle \prod_{i=1}^{i=q}(\tilde{\sigma}_1^{r_i} \tilde{\sigma}_1^{s_i} R^-(\tilde{\sigma}^{r_i}, \tilde{\sigma}^{s_i})) \right\rangle_-$$

$$= \prod_{i=1}^{i=q} \langle \tilde{\sigma}_1^{r_i} \tilde{\sigma}_1^{s_i} R^-(\tilde{\sigma}^{r_i}, \tilde{\sigma}^{s_i}) \rangle_-$$
$$= (\langle \tilde{\sigma}_1^1 \tilde{\sigma}_1^2 R^-(\tilde{\sigma}^1, \tilde{\sigma}^2) \rangle_-)^q.$$

La quantité $\langle f^- \prod_{i=1}^{i=q} R^-(\tilde{\sigma}^{r_i}, \tilde{\sigma}^{s_i}) \rangle_0$ ne dépend donc pas de la partition $\pi_q^*$; puisqu'il y a $\frac{(2q)!}{2^q}$ telles partitions ordonnées distinctes, en prenant l'espérance et en utilisant la seconde partie de la relation (9), on obtient l'égalité (36) annoncée. $\square$

4.2. *Développements de Taylor à l'ordre p et calcul des moments.* Suite à la Remarque 4.3, nous nous attacherons maintenant à appliquer la formule de Taylor pour calculer $\nu(f^- \tilde{\varepsilon}^{\otimes p})$ lorsque $f^- = \tilde{\sigma}_1^{\otimes p}$, en développant $\nu_t(f^- \tilde{\varepsilon}^{\otimes p})$ pour $t=1$ en $t=0$. Dans le cas $f^-$ quelconque, le Théorème 3.3 montre que le premier terme non nul du développement de Taylor est celui d'ordre $q$ lorsque $p = 2q$ ou $2q - 1$, le reste étant, d'après la relation (10) ou le Corollaire 2.4, un $O(q+1)\nu^{1/2}((f^-)^2)$. Nous montrerons que dans le cas particulier de notre $f^-$, on peut améliorer ce résultat, en développant à l'ordre $p$.

Plus précisément, nous montrerons le théorème suivant:

THÉORÈME 4.4. *Si $2l > p$,*
$$(37) \qquad \nu_0^{(l)}(\tilde{\varepsilon}^{\otimes p} \tilde{\sigma}_1^{\otimes p}) = O(p+1).$$

On en déduit le corollaire suivant:

COROLLAIRE 4.5. *On a*
$$(38) \quad \nu(\tilde{\varepsilon}^{\otimes p} \tilde{\sigma}_1^{\otimes p}) = \begin{cases} O(p+1), & si\ p = 2q-1, \\ 1/(q!)\nu_0^{(q)}(\tilde{\varepsilon}^{\otimes 2q} \tilde{\sigma}_1^{\otimes 2q}) + O(p+1), & si\ p = 2q. \end{cases}$$

DÉMONSTRATION DU COROLLAIRE 4.5. Il suffit de développer $\nu_t(\tilde{\varepsilon}^{\otimes p} \tilde{\sigma}_1^{\otimes p})$, pour $t=1$, par la formule de Taylor en $t=0$ jusqu'à l'ordre $p$, les dérivées en $t=0$ d'ordre $< q$ sont nulles par la première partie de la relation (9), les dérivées en $t=0$ d'ordre $> q (\geq q$ si $p = 2q-1)$ sont des $O(p+1)$ par le Théorème 4.4, la dérivée d'ordre $p+1$ en $t$ étant d'ordre $p+1$, par le Corollaire 2.4. $\square$

REMARQUE 4.6. On voit donc que si $p = 2q - 1$, le moment d'ordre p de la variable $N^{1/2} \tilde{\gamma}_{i,j}$ tend vers 0 si $N \to +\infty$, et si $p = 2q$, se comporte comme $\frac{N^q}{q!} \nu_0^{(q)}(\tilde{\varepsilon}^{\otimes 2q} \tilde{\sigma}_1^{\otimes 2q})$ si $N \to +\infty$. Nous entreprendrons l'étude de cette quantité dans la section suivante.

Afin de prouver ce théorème, énonçons tout d'abord un résultat préliminaire, que nous admettrons et dont la démonstration sera donnée dans l'Annxe B.

LEMME 4.7. *Soit $J_i$ une suite de $m$ doublets d'entiers de type $\{r_i, s_i\}$. Posons*

$$\hat{R}_i^- = \frac{1}{N} \sum_{j=1}^{j=N-2} \sigma_j^{r_i} \sigma_j^{s_i} - q_2. \tag{39}$$

*Alors on a*

$$\nu_-^{1/2}\left(\prod_{i=1}^{i=m} |\hat{R}_i^-|\right)^2 = O(m). \tag{40}$$

REMARQUE 4.8. Le symbole "$-$" dans $\hat{R}_i^-$ indique la dépendance en les $N-2$ premiers sites.

DÉMONSTRATION DU THÉORÈME 4.4. Plaçons nous dans le cas où $p = 2q$ ou $2q-1, 2l > p$ et posons $l = q + r$, avec $r \geq 1$ si $p$ est pair, $r \geq 0$ si $p$ est impair. On sait que, d'après le Corollaire 2.4, si $f^-$ est bornée, ce qui est le cas pour $f^- = \tilde{\sigma}_1^{\otimes p}$,

$$\nu_0^{(l)}(\tilde{\varepsilon}^{\otimes p} f^-) = O(q+r)\nu^{1/2}((f^-)^2) = O(p+1) \qquad \text{si } q + r \geq p + 1.$$

Le théorème est donc vrai pour de tels $r$.

Etudions alors cette expression lorsque $r \leq p - q$, c'est-à-dire $r \leq q - 1$ ou $q$, selon que $p$ est impair ou pair. On sait, d'après le Théorème 2.3 et la relation (71), que $\nu_0^{(q+r)}(\tilde{\varepsilon}^{\otimes p} f^-)$ est une combinaison linéaire finie de termes du type

$$\nu_0\left(\tilde{\varepsilon}^{\otimes p} \prod_{i=1}^{i=q+r} \varepsilon^{J_i}\right) \nu_0\left(f^- \prod_{i=1}^{i=q+r} (\dot{R}_{J_i}^-)\right).$$

Dans chacun de ces produits de deux termes, le premier est borné, car $|\tilde{\varepsilon}| \leq 2$. Nous allons étudier les termes $\nu_0(f^- \prod_{i=1}^{i=q+r}(\dot{R}_{J_i}^-))$ lorsque $f^- = \tilde{\sigma}_1^{\otimes p}$. Remarquons qu'alors $f^-$ est bornée.

On peut écrire, en utilisant la relation (71) et la définition de $\nu_-$ donnée dans l'Annexe A, la relation

$$\nu_0\left(\tilde{\sigma}_1^{\otimes p} \prod_{i=1}^{i=q+r} (\dot{R}_{J_i}^-)\right) \stackrel{\text{def}}{=} \nu_0(\tilde{\sigma}_1^{\otimes p} g^-) \stackrel{\text{def}}{=} \nu_0(\hat{h}) = \nu_-(\hat{h}).$$

Dans cette formule, $\hat{h}$ dépend uniquement des spins aux $N-1$ premiers sites. On peut alors montrer la proposition suivante:

PROPOSITION 4.9. *On a*

$$\nu_-(\hat{h}) = \nu_-(\tilde{\sigma}_i^{\otimes p} \times g^-) \qquad \forall i \leq N-1. \tag{41}$$

La démonstration se fait en utilisant la symétrie des $N-1$ premiers sites vis à vis de la mesure associée à $\nu_-$. □

Si l'on pose $\sigma_{N-1} = \bar{\varepsilon}$ et $\tilde{\sigma}_{N-1} = \hat{\varepsilon}$, et si de plus l'on remarque que lorsque $J_i = \{r_i, s_i\}$, en utilisant la définition de $\hat{R}_i^-$ donnée par la relation (39), on a

$$\nu_-(\hat{h}) = \nu_-(\hat{\varepsilon}^{\otimes p} \times g^-) = \nu_-\left[\hat{\varepsilon}^{\otimes p} \times \prod_{i=1}^{i=q+r}\left(\hat{R}_i^- + \frac{(\bar{\varepsilon})^{J_i}}{N}\right)\right]. \tag{42}$$

Si l'on développe le produit, on obtient

$$\nu_-(\hat{h}) = \sum_{k=0}^{k=q+r} \frac{1}{N^k} \sum_{\substack{B_1 \subset (1,2,\ldots,q+r) \\ |B_1|=k}} \nu_-\left(\hat{\varepsilon}^{\otimes p} \prod_{i \in B_1} (\bar{\varepsilon})^{J_i} \prod_{i \in B_1^c} \hat{R}_i^-\right). \tag{43}$$

Nous pouvons établir la proposition suivante, en utilisant notamment le Corollaire 3.5:

PROPOSITION 4.10. *Pour tout entier $k$ tel que $0 \leq k \leq q+r$, pour tout sous ensemble $B_1$ tel que $|B_1| = k$, on a*

$$\frac{1}{N^k}\nu_-\left(\hat{\varepsilon}^{\otimes p} \prod_{i \in B_1} (\bar{\varepsilon})^{J_i} \prod_{i \in B_1^c} \hat{R}_i^-\right) = O(p+1). \tag{44}$$

REMARQUE 4.11. Cette proposition une fois montrée, le Théorème 4.4 sera complètement démontré modulo la preuve du Lemme 4.7.

DÉMONSTRATION DE LA PROPOSITION 4.10. On sait, d'après le Lemme 4.7, où l'on pose $m = q+r-k$, que $\nu_-^{1/2}(\prod_{i \in B_1^c} |\hat{R}_i^-|)^2 = O(q+r-k)$, le terme correspondant à $k$ dans la somme ci dessus sera un $O(2k)O(q+r-k) = O(q+r+k) = O(p+1)$ si $k \geq p+1-q-r$.

Reste donc à étudier les cas des entiers $k$ vérifiant $k \leq p-q-r$. Nous appliquerons le Corollaire 3.5, utilisé pour $N-1$ sites, à

$$\nu_-\left(\hat{\varepsilon}^{\otimes p} \prod_{i \in B_1} (\bar{\varepsilon})^{J_i} \prod_{i \in B_1^c} \hat{R}_i^-\right) \stackrel{\text{def}}{=} \nu_-(u). \tag{45}$$

Si l'on pose $\prod_{i \in B_1} (\bar{\varepsilon})^{J_i} = (\bar{\varepsilon})^{B'}$, $B'$ vérifie alors les conditions du Corollaire 3.5, à savoir $|B'| \leq 2k < p$, et donc les dérivées d'ordre $< q-k$ de $(\nu_-)_t(u)$ en $t=0$ sont nulles, la dérivée d'ordre $q-k$ en $t=0$ sera majorée, puisque $\beta^- < \beta \leq \beta_0$, par

$$\frac{K}{(N-1)^{(q-k)/2}}\nu_-^{1/2}(u^2).$$

Comme $\frac{1}{N-1} < \frac{2}{N}$ et que, d'après le Lemme 4.7 (appliqué à $m = q-k+r$) et la définition de $u$,

$$\nu_-^{1/2}(u^2) = O(q+r-k),$$

on voit que cette dérivée est d'ordre $O(q-k)O(q+r-k) = O(2q+r-2k)$, le reste associé au développement de Taylor à l'ordre $q-k$ de $(\nu_-)_t(u)$ étant, par un raisonnement du même ordre et par le Corollaire 2.4, un

$$O(q-k+1)\nu_-^{1/2}\left(\prod_{i \in B_1^c} |\hat{R}_i^-|\right)^2 = O(2q+r+1-2k).$$

On en déduit que, par son développement de Taylor à l'ordre $q-k$, $\nu_-(u)$ est un $O(2q+r-2k)$. L'ordre de

$$\frac{1}{N^k} \sum_{\substack{B_1 \subset (1,2,\ldots,q+r) \\ |B_1|=k}} \nu_-\left(\hat{\varepsilon}^{\otimes p} \prod_{i \in B_1} (\bar{\varepsilon})^{J_i} \prod_{i \in B_1^c} \hat{R}_i^-\right)$$

vaut $O(2q+r)$ et donc chacun de ces termes est d'ordre $O(p+1)$. La Proposition 4.10 est complètement démontrée. $\square$

**5. Evaluation des moments d'ordre pair de la covariance symétrisée.** Certaines propositions de cette section seront démontrées dans l'Annexe B.

Pour évaluer, comme nous y invite la Remarque 4.6, le comportement asymptotique des moments d'ordre $2q$ de la variable $\tilde{\gamma}_{i,j}$, il suffit d'évaluer, d'après les relations (35) et (38), celui de $\nu_0^{(q)}(\tilde{\sigma}_1^{\otimes 2q}\tilde{\varepsilon}^{\otimes 2q})$, et donc, en raison de la relation (36), celui de la quantité

(46) $$E[\langle \tilde{\sigma}_1^1 \tilde{\sigma}_1^2 R^-(\tilde{\sigma}_1,\tilde{\sigma}_2)\rangle_-^q].$$

Nous parviendrons ainsi à notre résultat principal, qui s'énonce ainsi:

THÉORÈME 5.1. *Pour $\beta$ suffisamment petit,*

(47) $$\nu_0^{(q)}(\tilde{\sigma}_1^{\otimes 2q}\tilde{\varepsilon}^{\otimes 2q}) = (A_q')^2 \frac{(2q)!}{2^q}\left(\frac{4\beta^2}{N(1-\beta^2 A_1')}\right)^q + O(2q+1),$$

*et donc, d'après le Corollaire 4.5,*

(48) $$\nu(\tilde{\sigma}_1^{\otimes 2q}\tilde{\varepsilon}^{\otimes 2q}) = (A_q')^2 \frac{(2q)!}{q!2^q}\left(\frac{4\beta^2}{N(1-\beta^2 A_1')}\right)^q + O(2q+1).$$

On déduit immédiatemment du Théorème 5.1 et des équations (34) et (35) le résultat suivant:

COROLLAIRE 5.2. *Pour $\beta$ suffisamment petit, le moment d'ordre $2q$ de la variable $N^{1/2}\tilde{\gamma}_{i,j}$ est égal à*

$$(49) \qquad (A'_q)^2 \frac{(2q)!}{q! 2^q} \left( \frac{4\beta^2}{1-\beta^2 A'_1} \right)^q + O(1).$$

Pour montrer le Théorème 5.1, nous passerons par plusieurs étapes, que nous décrivons succinctement:

(50) Soit $\pi^*_{q,0}$ la partition canonique $\{1,2\}, \{3,4\}, \ldots, \{2q-1, 2q\}$.

Nous montrons tout d'abord que la quantité (46) s'écrit sous la forme

$$\nu_-(g_1^-), \qquad \text{où } g_1^- \stackrel{\text{def}}{=} \tilde{\sigma}_1^{\otimes 2q} R^-_{\pi^*_{q,0}},$$

$g_1^-$ dépend seulement des spins aux $N-1$ premiers sites.

Soit $g_1 \stackrel{\text{def}}{=} \tilde{\sigma}_1^{\otimes 2q} R_{\pi^*_{q,0}}$ un analogue de $g_1^-$ dépendant des spins à tous les $N$ sites. En utilisant notamment le Corollaire 3.4, nous obtenons une première évaluation de $\nu(g_1)$. Cette évaluation fait apparaitre des quantités

$$\nu_0(R^-_{\pi^*_{q-k,0}} R^-_{\hat{\pi}_{q-k}}) \stackrel{\text{def}}{=} \nu_-(g_2^-).$$

Comme çi dessus, $g_2^-$ dépend seulement des spins aux $N-1$ premiers sites. Soit $g_2$ un analogue de $g_2^-$ dépendant de tous les $N$ sites. Un TCL pour les recouvrements permet alors l'évaluation de $\nu(g_2)$ (la démonstration de la proposition correspondante est faite à la fin de la section). Cette évaluation permet alors celle de $\nu_-(g_2^-)$ (résultat montré dans l'Annexe B), ce qui permet d'obtenir un théorème évaluant $\nu(g_1)$ complétant la première évaluation déja obtenue. Une dernière proposition montre alors que $\nu(g_1)$ et $\nu_-(g_1^-)$ ont le même comportement asymptotique, ce qui entraine le Théorème 5.1.

Venons en au détail des énoncés et démonstrations:

Une premier calcul de la quantité (46) est donné par la proposition suivante:

PROPOSITION 5.3. *Soit $R'_{\pi^*_{q,0}}$ l'analogue de $R_{\pi^*_{q,0}}$ pour les $N-1$ premiers sites. Alors on a*

$$E[\langle \tilde{\sigma}_1^1 \tilde{\sigma}_1^2 R^-(\tilde{\sigma}_1, \tilde{\sigma}_2) \rangle^q_-] = E[\langle \tilde{\sigma}_1^{\otimes 2q} R^-_{\pi^*_{q,0}} \rangle_-]$$

$$(51) \qquad\qquad = \nu_-(\tilde{\sigma}_1^{\otimes 2q} R^-_{\pi^*_{q,0}})$$

$$\qquad\qquad = \left(\frac{N-1}{N}\right)^q \nu_-(\tilde{\sigma}_1^{\otimes 2q} R'_{\pi^*_{q,0}}).$$

DÉMONSTRATION. Les premières égalités de la relation (51) résultent du théorème de Fubini et de la définition de $\nu_-$, la dernière de la relation: Si

$$R'(\tilde\sigma^1,\tilde\sigma^2) = \frac{1}{N-1}\sum_{j=1}^{j=N-1}\tilde\sigma_j^1\tilde\sigma_j^2, \qquad R^-(\tilde\sigma^1,\tilde\sigma^2) = \frac{N-1}{N}R'(\tilde\sigma^1,\tilde\sigma^2).$$

On a alors

$$R^-_{\pi_{q,0}^*} = \prod_{i=1}^{i=q} R^-(\tilde\sigma^{2i-1},\tilde\sigma^{2i}) = \left(\frac{N-1}{N}\right)^q R'_{\pi_{q,0}^*}. \qquad \square$$

La Proposition 5.3 conduit donc à évaluer la quantité $\nu_-(\tilde\sigma_1^{\otimes 2q} R'_{\pi_{q,0}^*})$, dépendant des $N-1$ premiers sites.

Essayons d'évaluer plutot la quantité analogue construite sur $N$ sites, c'est à dire $\nu(\tilde\sigma_1^{\otimes 2q} R_{\pi_{q,0}^*})$. On obtient une première évaluation de cette quantité grâce à la proposition suivante:

PROPOSITION 5.4. *Posons*

(52) $$S_{0,k} = \sum_{\hat\pi_{q-k}} \nu_0(R^-_{\pi_{q-k,0}^*} R^-_{\hat\pi_{q-k}}).$$

*On a alors*

(53) $$\nu(\tilde\sigma_1^{\otimes 2q} R_{\pi_{q,0}^*}) = A'_q\left(\sum_{k=0}^{k=q}\left(\frac{4}{N}\right)^k C_q^k (\beta^2)^{q-k} S_{0,k}\right) + O(2q+1).$$

DÉMONSTRATION. Il est clair que, par symétrie des sites, en posant $\tilde\varepsilon^r = \tilde\sigma_N^r$, on a

$$\forall i \leq N \qquad \nu(\tilde\sigma_1^{\otimes 2q} R_{\pi_{q,0}^*}) = \nu(\tilde\sigma_i^{\otimes 2q} R_{\pi_{q,0}^*})$$
$$= \nu(\tilde\varepsilon^{\otimes 2q} R_{\pi_{q,0}^*}).$$

D'autre part,

$$\nu(\tilde\varepsilon^{\otimes 2q} R_{\pi_{q,0}^*}) = \nu\left(\tilde\varepsilon^{\otimes 2q} \prod_{i=1}^{i=q}\left(R^-(\tilde\sigma^{2i-1},\tilde\sigma^{2i}) + \frac{\tilde\varepsilon^{2i-1}\tilde\varepsilon^{2i}}{N}\right)\right)$$
$$= \sum_{k=0}^{k=q}\left(\frac{1}{N}\right)^k \times \sum_{\substack{B\subset(1,2\ldots,q)\\|B|=k}} \nu\left(\tilde\varepsilon^{\otimes 2q}\prod_{i\in B}(\tilde\varepsilon^{2i-1}\tilde\varepsilon^{2i})\right.$$

(54) $$\left.\times \prod_{i\notin B} R^-(\tilde\sigma^{2i-1},\tilde\sigma^{2i})\right),$$

$$= \sum_{k=0}^{k=q} \left(\frac{1}{N}\right)^k \times \sum_{\substack{B\subset(1,2,\ldots,q)\\|B|=k}} \nu\Bigg(\prod_{i\notin B}(\tilde{\varepsilon}^{2i-1}\tilde{\varepsilon}^{2i})\prod_{i\in B}(\tilde{\varepsilon}^{2i-1}\tilde{\varepsilon}^{2i})^2$$
$$\times \prod_{i\notin B} R^-(\tilde{\sigma}^{2i-1},\tilde{\sigma}^{2i})\Bigg).$$

Il est clair que, $k$ étant fixé, par invariance de $\nu$ par permutation des répliques, le $\nu(\cdot)$ associé à un $B$ donné tel que $|B| = k$ ne va pas dépendre de $B$; on peut donc réordonner les indices et supposer que

$$B^c = (1,2,\ldots,q-k), \qquad B = (q-k+1, q-k+2,\ldots,q).$$

Dans ce cas, en posant

(55) $$k' = q - k, \qquad g^- \stackrel{\text{def}}{=} \prod_{i=1}^{i=k'} R^-(\tilde{\sigma}^{2i-1},\tilde{\sigma}^{2i}) = R^-_{\pi^*_{k',0}},$$

on obtient

(56)
$$\nu\Bigg(\prod_{i\notin B}(\tilde{\varepsilon}^{2i-1}\tilde{\varepsilon}^{2i})\prod_{i\in B}(\tilde{\varepsilon}^{2i-1}\tilde{\varepsilon}^{2i})^2 \prod_{i\notin B} R^-(\tilde{\sigma}^{2i-1},\tilde{\sigma}^{2i})\Bigg)$$
$$= \nu\Bigg(\tilde{\varepsilon}^{\otimes 2k'}\prod_{i=2k'+1}^{i=2q}(\tilde{\varepsilon}^i)^2 g^-\Bigg).$$

En vertu du Corollaire 3.4 appliqué à $k'$ au lieu de $q$ et à

$$h_1(\eta_C) = \prod_{i=2k'+1}^{i=2q}(\tilde{\varepsilon}^i)^2, \qquad \text{où } C = (4k'+1, 4k'+2,\ldots,4q),$$

on obtient

(57)
$$\nu\Bigg(\tilde{\varepsilon}^{\otimes 2k'}\prod_{i=2k'+1}^{i=2q}(\tilde{\varepsilon}^i)^2 g^-\Bigg) = \beta^{2k'} E\Bigg[\bigg(\frac{1}{ch(Y)}\bigg)^{4k'}(\langle(\tilde{\varepsilon}^i)^2\rangle_0)^{2k}\Bigg]$$
$$\times \sum_{\hat{\pi}_{k'}} \nu_0(g^- R^-_{\hat{\pi}_{k'}})$$
$$+ O(k'+1)\nu^{1/2}((g^-)^2).$$

Le terme complémentaire de l'équation (57) est un $O(2k'+1)$, car on peut montrer, en utilisant la définition de $g^-$, les inégalités exponentielles (78) et les inégalités de Hölder que l'on a

$$\nu^{1/2}((g^-)^2) = O(q-k) = O(k').$$

Remarquons que la relation (52) entraine que

$$S_{0,k} = \sum_{\hat{\pi}_{k'}} \nu_0(g^- R^-_{\hat{\pi}_{k'}}).$$

Puisque

$$\langle (\tilde{\varepsilon}^i)^2 \rangle_0 = \langle (\varepsilon^{2i-1} - \varepsilon^{2i})^2 \rangle_0 = \langle 2(1 - \varepsilon^{2i-1}\varepsilon^{2i}) \rangle_0 = 2(1 - th^2(Y)) = \frac{2}{ch^2(Y)},$$

on a

$$E\left[\left(\frac{1}{ch(Y)}\right)^{4k'} (\langle (\tilde{\varepsilon}^i)^2 \rangle_0)^{2k}\right] = E\left[\left(\frac{1}{ch(Y)}\right)^{4k'} \left(\frac{2}{ch^2(Y)}\right)^{2k}\right]$$
$$= 4^k E\left(\frac{1}{ch^{4q}(Y)}\right) = 4^k A'_q.$$

Donc, d'après la relation (57), $\forall B \subset \{1, 2, \ldots, q\}$, tel que $|B| = k$, on a

(58)
$$\nu\left(\tilde{\varepsilon}^{\otimes 2q} \prod_{i \in B} \tilde{\varepsilon}^{2i-1}\tilde{\varepsilon}^{2i} \prod_{i \notin B} R^-(\tilde{\sigma}^{2i-1}, \tilde{\sigma}^{2i})\right)$$
$$= ((\beta^2)^{q-k} 4^k S_{0,k}) A'_q + O(2q - 2k + 1).$$

On a donc, en utilisant les relations (54) et (58),

$$\nu(\tilde{\sigma}_1^{\otimes 2q} R_{\pi^*_{q,0}}) = \left(\sum_{k=0}^{k=q} \left(\frac{4}{N}\right)^k C_q^k (\beta^2)^{q-k} S_{0,k}\right) A'_q + O(2q+1).$$

C'est l'équation (53) annoncée dans l'énoncé de la proposition. □

La Proposition 5.4 conduit donc à devoir évaluer les sommes $S_{0,k}$ définies dans la relation (52), ou plutôt, pour toute partition $\hat{\pi}_{q-k}$, leurs termes $\nu_0(R^-_{\pi^*_{q-k,0}} R^-_{\hat{\pi}_{q-k}})$. Ces termes dépendent des $N-1$ premiers sites, essayons d'évaluer leurs analogues pour $N$ sites, soit les nombres $\nu(R_{\pi^*_{q-k,0}} R_{\hat{\pi}_{q-k}})$. Leur comportement va être décrit par la proposition qui suit:

PROPOSITION 5.5. *On a*

(59)
$$\nu(R_{\hat{\pi}_{q-k}} R_{\pi^*_{q-k,0}})$$
$$= \begin{cases} O(2q - 2k + 1), & \text{si } \hat{\pi}_{q-k} \neq \pi^*_{q-k,0}, \\ \left(\frac{4A'_1}{N(1 - \beta^2 A'_1)}\right)^{q-k} + O(2(q-k) + 1), & \text{si } \hat{\pi}_{q-k} = \pi^*_{q-k,0}. \end{cases}$$

Cette proposition, dont la démonstration s'appuie sur la démonstration de théorèmes centraux limites pour les recouvrements $R$, sera provisoirement admise. Elle sera démontrée à la fin de cette sous section.

On montrera également, cette fois dans l'Annexe B, la proposition suivante:

PROPOSITION 5.6. *Si l'on pose $k' = q - k$, on a*

$$\nu_0(R^-_{\pi^*_{k',0}} R^-_{\hat{\pi}_{k'}}) = \nu(R_{\pi^*_{k',0}} R_{\hat{\pi}_{k'}}) + O(2k' + 1). \tag{60}$$

Ces deux propositions, une fois démontrées, permettent d'obtenir une évaluation asymptotique précise de $\nu(\tilde{\sigma}_1^{\otimes 2q} R_{\pi^*_{q,0}})$, donnée par le théorème suivant:

THÉORÈME 5.7. *On a*

$$\nu(\tilde{\sigma}_1^{\otimes 2q} R_{\pi^*_{q,0}}) = \left(\frac{4}{N(1 - \beta^2 A'_1)}\right)^q A'_q + O(2q + 1). \tag{61}$$

DÉMONSTRATION. Commençons par compléter la Proposition 5.4, en calculant, par application des Propositions 5.5 et 5.6, les termes des sommes $S_{0,k}$ introduites dans l'équation (52). On obtient

$$\begin{aligned}
\nu_0(R^-_{\pi^*_{k',0}} R^-_{\hat{\pi}_{k'}}) &= \nu(R_{\pi^*_{k',0}} R_{\hat{\pi}_{k'}}) + O(2(q-k) + 1) \\
&= \begin{cases} O(2k' + 1), & \text{si } \hat{\pi}_{k'} \neq \pi^*_{k',0}, \\ \left(\dfrac{4A'_1}{N(1 - \beta^2 A'_1)}\right)^{k'} + O(2k' + 1), & \text{si } \hat{\pi}_{k'} = \pi^*_{k',0}. \end{cases}
\end{aligned} \tag{62}$$

Comme

$$S_{0,k} = \sum_{\hat{\pi}_{k'}} \nu_0(R^-_{\pi^*_{k',0}} R^-_{\hat{\pi}_{k'}}),$$

on a

$$S_{0,k} = \left(\frac{4A'_1}{N(1 - \beta^2 A'_1)}\right)^{q-k} + O(2(q - k) + 1).$$

Il s'ensuit que

$$\begin{aligned}
\nu(\tilde{\sigma}_1^{\otimes 2q} R_{\pi^*_{q,0}}) &= A'_q \sum_{k=0}^{k=q} C^k_q \left(\frac{4}{N}\right)^k (\beta^2)^{q-k} \left(\frac{4A'_1}{N(1 - \beta^2 A'_1)}\right)^{q-k} + O(2q + 1) \\
&= A'_q \left(\frac{4}{N}\right)^q \left(\sum_{k=0}^{k=q} C^k_q \left(\frac{\beta^2 A'_1}{1 - \beta^2 A'_1}\right)^{q-k}\right) + O(2q + 1)
\end{aligned}$$

$$= A'_q \left(\frac{4}{N}\right)^q \left(1 + \frac{\beta^2 A'_1}{1 - \beta^2 A'_1}\right)^q + O(2q+1)$$

$$= \left(\frac{4}{N(1 - \beta^2 A'_1)}\right)^q A'_q + O(2q+1).$$

C'est la formule (61) annoncée. □

Nous pouvons maintenant démontrer notre résultat principal:

DÉMONSTRATION DU THÉORÈME 5.1. On sait, d'après la relation (36) et la Proposition 5.3, que

$$\nu_0^{(q)}(\tilde{\sigma}_1^{\otimes 2q} \tilde{\varepsilon}^{\otimes 2q}) = \beta^{2q} A'_q \frac{(2q)!}{2^q} E[(\langle \tilde{\sigma}_1^1 \tilde{\sigma}_1^2 R^-(\tilde{\sigma}^1, \tilde{\sigma}^2)\rangle_-)^q]$$

$$= \beta^{2q} A'_q \frac{(2q)!}{2^q} \nu_-(\tilde{\sigma}_1^{\otimes 2q} R^-_{\pi^*_{q,0}}).$$

Le théorème sera démontré dès que l'on aura montré le résultat suivant:

PROPOSITION 5.8. *Les quantités* $\nu_-(\tilde{\sigma}_1^{\otimes 2q} R^-_{\pi^*_{q,0}})$ *et* $\nu(\tilde{\sigma}_1^{\otimes 2q} R_{\pi^*_{q,0}})$ *ont la même évaluation asymptotique, donnée par le Théorème* 5.7, *et l'on a donc*

$$(63) \qquad \nu_-(\tilde{\sigma}_1^{\otimes 2q} R^-_{\pi^*_{q,0}}) = \left(\frac{4}{N(1-\beta^2 A'_1)}\right)^q A'_q + O(2q+1).$$

Ce résultat sera montré dans l'Annexe B.

DÉMONSTRATION DE LA PROPOSITION 5.5. Avant d'aborder la démonstration proprement dite, montrons deux lemmes préliminaires. Utilisons les notations (79) développées dans l'Annexe A à propos de la démonstration des TCL sur les recouvrements.

LEMME 5.9. *Si* $\tilde{\sigma}^r = \sigma^{2r-1} - \sigma^{2r}$ *et* $\tilde{\sigma}^s = \sigma^{2s-1} - \sigma^{2s}$,

$$(64) \qquad R(\tilde{\sigma}^r, \tilde{\sigma}^s) = T_{2r-1, 2s-1} + T_{2r, 2s} - T_{2r-1, 2s} - T_{2r, 2s-1}.$$

DÉMONSTRATION. On a

$$R(\tilde{\sigma}^r, \tilde{\sigma}^s) = R_{2r-1, 2s-1} - q_2 + R_{2r, 2s} - q_2 - (R_{2r-1, 2s} - q_2) - (R_{2r, 2s-1} - q_2).$$

On sait que

$$R_{i,j} - q_2 = T_{i,j} + T_i + T_j + T.$$

On a bien le résultat désiré. □

REMARQUE 5.10. L'équation (64) s'écrit aussi

$$R(\tilde{\sigma}^r, \tilde{\sigma}^s) = \sum_{\substack{l \in \{2r-1, 2r\} \\ l' \in \{2s-1, 2s\}}} (-1)^{l+l'} T_{l,l'}. \tag{65}$$

Nous aurons également besoin d'utiliser la variante du Théorème A.5 suivante:

LEMME 5.11. *Plaçons nous dans les conditions du Théorème A.5. Soit $\Delta$ un ensemble de doublets $\{l, l'\}$, où $1 \leq l < l' \leq n$. Soit $k(l, l')$ un entier $\geq 1$ associé à chacun de ces doublets, $k^+$ la somme de ces entiers. Soit $a(r)$ le moment d'ordre $r$ d'une gaussienne centrée réduite. Soit*

$$A^2 = \frac{A_1'}{N(1 - \beta^2 A_1')} \quad et \quad A_1' = E\left(\frac{1}{ch^4(Y)}\right).$$

*Si $\beta \leq \beta_0$, on a*

$$\nu\left(\prod_{\{l,l'\} \in \Delta} T_{l,l'}^{k_{l,l'}}\right) = \prod_{\{l,l'\} \in \Delta} (a(k(l,l'))) A^{k^+} + O(k^+ + 1). \tag{66}$$

*Si, de plus, les entiers $k(l, l')$ valent 1 ou 2, le terme principal du membre de droite de l'équation (66) vaut 0, sauf si $\forall \{l, l'\} \in \Delta, k(l, l') = 2$. On a dans ce cas*

$$\nu\left(\prod_{\{l,l'\} \in \Delta} T_{l,l'}^2\right) = A^{k^+} + O(k^+ + 1). \tag{67}$$

DÉMONSTRATION. La première égalité (66) est une conséquence de la relation (81) de l'Annexe A, dans la seconde partie, on utilise la propriété: $a(1) = 0$, $a(2) = 1$. □

Passons maintenant à la démonstration proprement dite de la Proposition 5.5. Soit $\{r_i, s_i\}$, $i \in (1, 2, \ldots, q-k)$ la partition ordonnée naturelle associée à $\hat{\pi}_{q-k}$. On voit que, en utilisant l'équation (65), le produit

$$R_{\hat{\pi}_{q-k}} R_{\pi_{q-k,0}^*} = \prod_{i=1}^{i=q-k} R(\tilde{\sigma}^{2i-1}, \tilde{\sigma}^{2i}) R(\tilde{\sigma}^{r_i}, \tilde{\sigma}^{s_i})$$

est la somme de $4^{2(q-k)}$ expressions, qui sont chacune un produit de $2(q-k)$ termes du type $(-1)^{l+l'} T_{l,l'}$, non necessairement distincts.

Toutefois, les termes $(-1)^{l+l'} T_{l,l'}$ associés à chaque partition $\hat{\pi}_{q-k}$ ou $\pi_{q-k,0}^*$ qui apparaissent dans une expression donnée sont nécessairement distincts, alors que le même terme peut venir de l'une ou l'autre des deux partitions et donc, pour les doublets distincts $\{l, l'\}$ apparaissant dans l'expression, $k(l,l')$ désignant le nombre de fois que le doublet $\{l,l'\}$ apparait dans le produit, on a $1 \leq k(l,l') \leq 2$ et $k^+ = 2(q-k)$. On a $\forall \{l,l'\}, k(l,l') = 2$ si et seulement si les partitions $\hat{\pi}_{q-k}$ et $\pi^*_{q-k,0}$ coincident et si de plus pour former l'expression on prend deux fois le même terme $(-1)^{l+l'}T_{l,l'}$ pour chaque doublet $\{2i-1,2i\}$ de la partition $\pi^*_{q-k,0}$. Le Lemme 5.11 entraine la première partie de la relation (59). Quant à la seconde partie, on voit, par l'équation (67), que chaque choix (répété une fois) d'un des quatre termes $(-1)^{l+l'}T_{l,l'}$ associés à chaque doublet $\{2i-1,2i\}$ de la partition $\pi^*_{q-k,0}$ donne, pour l'expression produit de tels termes, notés alors $t_i$

$$\nu\left(\prod_{i=1}^{i=q-k}(t_i^2)\right) = A^{2(q-k)} + O(2(q-k)+1).$$

Et donc, en tenant compte des quatre choix possibles pour chaque doublet:

$$\nu\left(\prod_{i=1}^{i=q-k} R(\tilde{\sigma}^{2i-1},\tilde{\sigma}^{2i})^2\right) = (4A^2)^{q-k} + O(2(q-k)+1).$$

C'est la seconde partie de la relation (59). □
□

**6. Loi limite de la covariance des spins en deux sites.** Dans cette section, nous donnerons la preuve du Théorème 1.1. Plaçons nous dans les conditions de l'énoncé de ce théorème et utilisons les notations qui y sont introduites. Nous prendrons $\beta$ suffisamment petit pour que les Théorèmes 5.1 et 4.4 soient vérifiés.

Il est clair que $(A'_q)$ est le moment d'ordre $2q$ de la variable $U_1$ ou $U_2$, et donc $(A'_q)^2$ est le moment d'ordre $2q$ du produit $W = U_1U_2$.

D'autre part, le moment d'ordre $2q$ d'une variable gaussienne centrée réduite $z$ est $\frac{2q!}{2^q q!}$. Le Corollaire 5.2 peut donc se lire ainsi: $z$ étant indépendant de $U_1$ et $U_2$, le moment d'ordre $2q$ de la variable $N^{1/2}\tilde{\gamma}_{i,j}$ est égal à celui de la variable $\frac{2\beta}{\sqrt{1-\beta^2 A'_1}}zU_1U_2 + O(1)$.

En utilisant la relation (34) $\gamma_{i,j} = \frac{1}{2}\tilde{\gamma}_{i,j}$, et pour le cas des moments d'ordre impair le Corollaire 4.5 et la Remarque 4.6, on obtient bien le résultat annoncé.

## ANNEXE A

Nous décrirons dans cette section quelques notions développées par Talagrand dans son livre [10] et fréquemment utilisées içi.

### A.1. Méthode de la cavité et applications.

*Formule de la cavité.* Soit $\sigma$ une $N$-configuration. On peut écrire

$$\sigma = (\rho, \varepsilon) \qquad \text{où } \rho = (\sigma_1, \ldots, \sigma_i, \ldots, \sigma_{N-1}) \text{ et } \varepsilon = \sigma_N.$$

Soit

$$\beta^- = \beta\sqrt{\frac{N-1}{N}}. \tag{68}$$

On a

$$-H_N(\sigma, \beta, h) = -H_{N-1}(\rho, \beta^-, h) + \varepsilon(g(\rho) + h),$$

$$\text{où } g(\rho) = \frac{\beta}{\sqrt{N}} \sum_{i \leq N-1} g_{i,N} \sigma_i.$$

NOTATIONS A.1. Nous introduirons les notations suivantes:

- Le variable $G_{N-1}^-$ est la mesure de Gibbs associée à $-H_{N-1}(\rho, \beta^-, h)$.
- Le terme $\langle f^- \rangle_-$ est l'intégrale de $f^-$ par rapport à la mesure de Gibbs produit $(G_{N-1}^-)^{\otimes k}$, $f^-$ étant une fonction de $k$ configurations dépendant seulement des $N-1$ premiers sites.
- Le variable $\nu_-(f^-)$ est l'espérance de $\langle f^- \rangle_-$.

Soit $Z = \langle ch(g(\rho) + h) \rangle_-$. Si $f$ est une fonction de $k$ configurations, il existe une fonction $\hat{f}^-$ liée à $f$ ne dépendant pas des spins au site $N$ telle que $\langle f \rangle = \langle \hat{f}^- \rangle_- / Z^k$; c'est la formule de la cavité, qui consiste donc à pouvoir passer de $N-1$ sites à $N$ sites (voir [10], 2.153), ce qui a ensuite conduit à élaborer la méthode du "smart path."

*Construction des hamiltoniens du "smart path."* Posons $\forall t \in [0, 1]$ (intervalle de réels),

$$g_t(\rho) = \sqrt{t} g(\rho) + \sqrt{1-t} \beta z \sqrt{q_2}.$$

Les hamiltoniens de la méthode du smart path sont définis par la formule

$$-H_{N,t}(\sigma, \beta, h) = -H_{N-1}(\rho, \beta^-, h) + \varepsilon(g_t(\rho) + h).$$

Dans cette formule, $q_2$ vérifie la relation (3).

Chacun de ces hamiltoniens engendre, nous l'avons vu, une mesure de Gibbs, une intégrale et son espérance que nous notons respectivement, $G_{N,t}$, $\langle f \rangle_t$ et $\nu_t(f) = E(\langle f \rangle_t)$, ces deux dernières invariantes par permutation des répliques. On a alors $\langle f \rangle_1 = \langle f \rangle$ et $\nu_1(f) = \nu(f)$.

*Étude de $G_{N,0}$.*  On a la relation, en utilisant l'équation (3):

$$-H_{N,0}(\sigma,\beta,h) = -H_{N-1}(\rho,\beta^-,h) + \varepsilon Y.$$

On en déduit que le spin au site $N$ est indépendant des spins aux $N-1$ premiers sites, pour la mesure de Gibbs $G_{N,0}$ comme pour son espérance. Les spins aux $N-1$ premiers sites ont alors la loi $G_{N-1}^-$. Le spin au site $N$ vérifie la relation

(69) $$\langle \varepsilon \rangle_0 = th(Y).$$

Ces propriétés s'étendent au cas de plusieurs répliques. On montre que si la fonction de $k$ répliques $f = f^- h_0$, où $f^-$ dépend seulement des spins de ces répliques aux $N-1$ premiers sites et $h_0$ de leurs spins au site $N$, on obtient les relations

(70) $$\langle f \rangle_0 = \langle f^- \rangle_0 \langle h_0 \rangle_0 = \langle f^- \rangle_- \langle h_0 \rangle_0,$$

(71) $$\nu_0(f) = \nu_0(f^-)\nu_0(h_0) = \nu_-(f^-)\nu_0(h_0).$$

*Majoration de $\nu_t(f)$* (*voir* [10], *Proposition* 2.4.6).  En utilisant l'expression (6) donnant la dérivée de $\nu_t(f)$, on obtient la proposition suivante:

PROPOSITION A.2.  *Soit $f$ une fonction numérique de $n$ répliques, à valeurs positives. On a l'inégalité*

(72) $$\nu_t(f) \leq \exp(4n^2\beta^2)\nu(f) \leq K\nu(f), \qquad si\ \beta \leq \beta_0;$$

*$K$ étant une constante ne dépendant que de $n$.*

Dans le deux sous sections qui suivent, nous utiliserons les notations $O(k)$ définies dans l'Introduction (voir Définition 1.3).

**A.2. Inégalités exponentielles.**  Dans [10] (Section 2.5), est obtenu le résultat suivant:

PROPOSITION A.3.  *Il existe $\beta_0 > 0$ et une constante $L > 4$ telle que, si $\beta < \beta_0$, on a*

(73) $$\forall k \geq 0 \qquad \nu[(R_{1,2}-q_2)^{2k}] \leq (Lk/N)^k.$$

*Il en résulte, avec nos notations, que*

(74) $$\nu[(R_{1,2}-q_2)^{2k}] = O(2k).$$

REMARQUE A.4.  Cette proprieté implique l'existence d'une constante $L'$ telle que $\nu(\exp(N/L'(R_{1,2}-q_2)^2)) \leq L'$, d'où l'origine du nom de ces inégalités.

INDICATIONS DE DÉMONSTRATION. La propriété est immédiate pour $k \geq N$, car $(R_{1,2}-q_2)^{2k} < 4^k$. Pour $k < N$, la démonstration se fait par récurrence sur $k$.

Pour $k=1$, elle résulte, par symétrie des sites, de l'égalité

$$\nu(R_{1,2}-q_2)^2 = \nu(\varepsilon^1\varepsilon^2 - q_2)(R_{1,2}-q_2) \stackrel{\text{def}}{=} \nu(f).$$

On trouve $\nu_0(f) = (1-q_2^2)/N$, et l'on utilise le théorème des accroissements finis pour majorer $|\nu(f) - \nu_0(f)|$, ce qui conduit à l'inégalité désirée pour $\beta$ suffisamment petit. On procède de manière analogue pour montrer la récurrence. □

A l'aide des inégalités de Hölder, on en déduit également les inégalités

(75) $$\nu(|R_{1,2}-q_2|^j) \leq (L(j+1)/2N)^{j/2}.$$

On a donc là aussi l'inégalité

(76) $$\nu(|R_{1,2}-q_2|)^j = O(j).$$

Enfin, dans le courant de la démonstration, il est montré l'inégalité

(77) $$\nu[(R_{1,2}^- - q_2)^{2k}] \leq 3(L(k+1)/N)^k = O(2k).$$

Soient $(\sigma^1, \sigma^2, \sigma^3, \sigma^4)$ quatre répliques, $\tilde{\sigma}^1 = \sigma^1 - \sigma^2, \tilde{\sigma}^2 = \sigma^3 - \sigma^4$. Posons

$$R(\tilde{\sigma}^1, \tilde{\sigma}^2) = 1/N \left( \sum_{i=1}^{i=N} \tilde{\sigma}_i^1 \tilde{\sigma}_i^2 \right)$$

et de même

$$R^-(\tilde{\sigma}^1, \tilde{\sigma}^2) = 1/N \left( \sum_{i=1}^{i=N-1} \tilde{\sigma}_i^1 \tilde{\sigma}_i^2 \right).$$

On a les relations

$$R^-(\tilde{\sigma}^1, \tilde{\sigma}^2) = (R^-(\sigma^1, \sigma^3) - q_2) + (R^-(\sigma^2, \sigma^4) - q_2)$$
$$- (R^-(\sigma^1, \sigma^4) - q_2) - (R^-(\sigma^2, \sigma^3) - q_2).$$

Il en résulte, en utilisant les propriétés des normes usuelles sur un espace $L^{2k}(\nu)$ et l'inégalité (77), que

(78) $$\nu(R^-(\tilde{\sigma}^1, \tilde{\sigma}^2)^{2k}) = O(2k).$$

**A.3. Théorèmes centraux limites pour les recouvrements (voir [10], Sections 2.6 et 2.7.).** Soit $(\sigma^1, \ldots, \sigma^n)$ une suite de $n$ répliques. Considérons les $C_n^2 \stackrel{\text{def}}{=} m$ recouvrements construits sur ces $n$ répliques, notés

$$R_{l,l'} = R(\sigma^l, \sigma^{l'}), \qquad \text{où } 1 \leq l < l' \leq n.$$

On peut alors définir les recouvrements centrés $R_{l,l'} - \langle R_{l,l'} \rangle$. Soit $X_{l,l'} = N^{1/2}(R_{l,l'} - \langle R_{l,l'} \rangle)$. Soit $X$ le vecteur dont les $m$ composantes sont les variables $X_{l,l'}$. Considérons les moments de tous ordres, pour la mesure de Gibbs $G_N^{\otimes n}$, du vecteur aléatoire $X$ sur $(\Sigma_N^{\otimes n}, G_N^{\otimes n})$, c'est-à-dire, pour toute suite de $m$ entiers $(k_{l,l'})$ les quantités

$$M_{(k_{l,l'})} \stackrel{\text{def}}{=} \left\langle \prod_{1 \leq l < l' \leq n} X_{l,l'}^{k_{l,l'}} \right\rangle.$$

Ces moments sont eux mêmes des variables aléatoires relativement au désordre du modèle.

Talagrand montre, dans la Section 2.7 de son livre [10], pour $\beta$ suffisamment petit, l'existence d'un vecteur gaussien centré de dimension $m$, dont il peut donner la matrice de covariance, et dont les moments d'ordre $(k_{l,l'})$ sont la limite dans $L^2$ des moments $M_{(k_{l,l'})}$. Il montre également la normalité asymptotique, au sens de la convergence des moments, de la variable $N^{1/2}(\langle R_{l,l'} \rangle - q_2)$ et donne la variance de la loi limite.

La démonstration de ces propriétés repose sur une décomposition des recouvrements "recentrés" $R_{l,l'} - q_2$. Rappelons quelques notions utilisées dans ce contexte (voir [10], formules 2.265 et 2.266): nous noterons $b_j = \langle \sigma_j \rangle, b$ le vecteur de $\mathbf{R}^N$ de composantes $b_j$; $\dot{\sigma}_j = \sigma_j - b_j, \dot{\sigma}$ le vecteur de $\mathbf{R}^N$ de composantes $\dot{\sigma}_j$; $R(u,v) = \frac{1}{N} \sum_{j=1}^{j=N} u_j v_j$, extension à $\mathbf{R}^N$ des recouvrements.

Soient $\sigma^l$ et $\sigma^{l'}$ deux $N$ configurations. On note

(79) $\qquad T_{l,l'} = R(\dot{\sigma}^l, \dot{\sigma}^{l'}), \qquad T_l = R(\dot{\sigma}^l, b), \qquad T = R(b,b) - q_2.$

On a

(80) $\qquad R_{l,l'} - q_2 = R(\sigma^l, \sigma^{l'}) - q_2 = T_{l,l'} + T_l + T_{l'} + T.$

La démonstration des TCL repose sur l'étude des moments de tous ordre du vecteur formé par toutes les variables introduites dans les relations (79), pour la mesure $\nu$, considérée içi comme espérance de la mesure de Gibbs produit $G_N^{\otimes n}$. On a en particulier le théorème suivant, qui est le Théorème 2.7.1 du livre [10]:

THÉORÈME A.5. *Soit $k(l,l')$ un entier $\geq 0$ associé à chacun des doublets $(l,l')$, avec $1 \leq l < l' \leq n$, soit $k^+$ la somme de ces entiers. Soit $a(r)$ le moment d'ordre $r$ d'une gaussienne centrée réduite. Soit*

$$A^2 = \frac{A_1'}{N(1 - \beta^2 A_1')} \qquad et \quad A_1' = E\left(\frac{1}{ch^4(Y)}\right).$$

*Si* $\beta \leq \beta_0$, *on a*

$$\nu\left(\prod_{(l,l')} T_{l,l'}^{k_{l,l'}}\right) = \prod_{(l,l')} (a(k(l,l')))A^{k^+} + O(k^+ + 1). \tag{81}$$

## ANNEXE B

DÉMONSTRATION DU LEMME 4.7. Il suffit de montrer, en utilisant les inégalités de Hölder comme dans la démonstration du Corollaire 2.4 que

$$\nu_-((\hat{R}_i^-)^{2m}) = O(2m).$$

On sait que, d'après la relation (72) de l'Annexe A, il existe une constante $K$ telle que

$$\nu_-((\hat{R}_i^-)^{2m}) = \nu_0((\hat{R}_i^-)^{2m}) \leq K\nu((\hat{R}_i^-)^{2m}),$$

$$(\hat{R}_i^-)^{2m} = \left(\dot{R}_{J_i}^- - \frac{\sigma_{N-1}^{r_i}\sigma_{N-1}^{s_i}}{N}\right)^{2m} \leq \sum_{l=0}^{l=2m} C_{2m}^l \frac{(|(\dot{R}_{J_i}^-)|)^l}{N^{2m-l}}.$$

En utilisant les inégalités exponentielles (77) appliquées aux quantités $\nu(\dot{R}_{J_i}^-)^{2l}$, on obtient

$$\nu_-((\hat{R}_i^-)^{2m}) \leq \sum_{l=0}^{l=2m} C_{2m}^l \frac{\nu(|\dot{R}_{J_i}^-|)^l}{N^{2m-l}}$$

$$\leq \sum_{l=0}^{l=2m} O(2(2m-l))O(l) = \sum_{l=0}^{l=2m} O(4m-l) = O(2m).$$

C'est la relation cherchée. □

DÉMONSTRATION DE LA PROPOSITION 5.6. Soit $\{r_i, s_i\}, i \in (1, 2, \ldots, k')$, la suite de doublets constituant la partition ordonnée naturelle (voir Remarque 3.2) associée à $\hat{\pi}_{k'}$:

$$\begin{aligned}
\nu(R_{\pi_{k',0}^*} R_{\hat{\pi}_{k'}}) &= \nu\left(\prod_{i=1}^{i=k'} \left(R^-(\tilde{\sigma}^{2i-1}, \tilde{\sigma}^{2i}) + \frac{\tilde{\varepsilon}^{2i-1}\tilde{\varepsilon}^{2i}}{N}\right)\right. \\
&\qquad\qquad \left. \times \left(R^-(\tilde{\sigma}^{r_i}, \tilde{\sigma}^{s_i}) + \frac{\tilde{\varepsilon}^{r_i}\tilde{\varepsilon}^{s_i}}{N}\right)\right) \\
&\stackrel{\text{def}}{=} \nu\left(\prod_{i=1}^{i=k'} (R^-(\tilde{\sigma}^{2i-1}, \tilde{\sigma}^{2i})R^-(\tilde{\sigma}^{r_i}, \tilde{\sigma}^{s_i}))\right) + \hat{T}_1 \\
&= \nu(R_{\pi_{k',0}^*}^- R_{\hat{\pi}_{k'}}^-) + \hat{T}_1 \stackrel{\text{def}}{=} \nu(\hat{f}) + \hat{T}_1.
\end{aligned} \tag{82}$$

Posons

(83) $$J_i \stackrel{\text{def}}{=} \begin{cases} \{2i-1, 2i\}, & \text{si } i \leq k', \\ \{r_i, s_i\}, & \text{si } k' < i \leq 2k', \end{cases}$$

(84) $$\tilde{\varepsilon}^J = \tilde{\varepsilon}^r \tilde{\varepsilon}^s \quad \text{si } J = \{r, s\},$$

(85) $$R_i^- \stackrel{\text{def}}{=} \begin{cases} R^-(\tilde{\sigma}^{2i-1}, \tilde{\sigma}^{2i}), & \text{si } i \leq k', \\ R^-(\tilde{\sigma}^{r_i}, \tilde{\sigma}^{s_i}), & \text{si } k' < i \leq 2k'. \end{cases}$$

L'expression $\hat{T}_1$ figurant dans (82) est alors une somme de termes du type

(86) $$\frac{1}{N^{|B|}} \nu\left(\prod_{i \in B} \tilde{\varepsilon}^{J_i} \prod_{i \notin B} R_i^-\right), \qquad \text{où } B \subset (1, 2, \ldots, 2k') \text{ et } |B| = l \geq 1.$$

Il y a $2k' - l$ termes du type $R_i^-$, et l'on peut montrer, en utilisant les inégalités exponentielles (78) de l'Annexe A que

$$\forall l \qquad \nu^{1/2}\left(\prod_{i \notin B}(R_i^-)\right)^2 = O(2k' - l).$$

Il s'ensuit que puisque ici $l \geq 1$, chaque terme $\frac{1}{N^{|B|}} \nu(\prod_{i \in B}(\tilde{\varepsilon}^{J_i}) \prod_{i \notin B} R_i^-)$ de l'équation (86) est un $O(2k' + l) = O(2k' + 1)$, et donc leur somme $\hat{T}_1 = O(2k' + 1)$.

On a, en utilisant [10], Proposition 2.5.3, et l'équation (82),

$$\nu(R_{\pi_{k',0}^*} R_{\hat{\pi}_{k'}}) = \nu(\hat{f}) + O(2k' + 1) = \nu_0(\hat{f}) + O(2k' + 1) + O(1)\nu^{1/2}(\hat{f}^2).$$

Par les inégalités exponentielles (78) et les inégalités de Hölder, on montre que $\nu^{1/2}(\hat{f}^2) = O(2k')$ et donc

$$\nu(R_{\pi_{k',0}^*} R_{\hat{\pi}_{k'}}) = \nu_0(\hat{f}) + O(2k' + 1).$$

C'est l'équation (60). □

DÉMONSTRATION DE LA PROPOSITION 5.8.    Nous devons donc évaluer la quantité

(87) $$\nu_-(\tilde{\sigma}_1^{\otimes 2q} R_{\pi_{q,0}^*}^-).$$

Rappelons que $\beta_-^2 = (\frac{N-1}{N})\beta^2$. Nous utiliserons également les notations suivantes:

$$Y^- = \beta_- \sqrt{q_2^-} z + h, \qquad \text{où } E[th^2(Y^-)] = q_2^-,$$
$$(A_q')^- = E[(1 - th^2(Y^-))^{2q}].$$

En utilisant la Proposition 5.3 puis le Théorème 5.7 avec $N-1$ remplaçant $N$ et $\beta_-$ remplaçant $\beta$, on voit que la quantité (87) vaut

$$\left(\frac{N-1}{N}\right)^q \nu_-(\tilde{\sigma}_1^{\otimes 2q} R'_{\pi_{q,0}^*})$$
$$(88) \qquad = \left(\frac{N-1}{N}\right)^q \left(\left(\frac{4}{(N-1)(1-\beta_-^2(A'_1)^-)}\right)^q (A'_q)^- + O^-(2q+1)\right).$$

La quantité $O^-(2q+1)$ est bornée par $\frac{K}{(N-1)^{q+1/2}}$ où $K$ est une constante ne dépendant pas de $\beta_-$ si $\beta_- \leq \beta_0$; puisque $\beta_- < \beta$, il suffit pour cela, ce que nous supposons, que l'on aie $\beta \leq \beta_0$. On a donc $O^-(2q+1) = O(2q+1)$.

On sait ([10], Lemme 2.4.15) que $q_2 - q_2^- = O(2)$ et donc les quantités nouvelles introduites $Y^-$ et $(A'_q)^-$ ont des limites si $N \to +\infty$, obtenues en supprimant les indices "$-$" dans leur definition.

Le terme principal du membre de droite de l'égalité (88) s'écrit

$$\left(\frac{1}{N}\right)^q \left(\frac{4}{1-\beta_-^2(A'_1)^-}\right)^q (A'_q)^-.$$

On veut montrer que ce terme s'écrit en fait

$$\left(\frac{1}{N}\right)^q \left(\frac{4}{1-\beta^2 A'_1}\right)^q A'_q + O(2q+1),$$

ce qui montrera le lemme. Il suffit de prouver que

$$\left(\frac{4}{1-\beta_-^2(A'_1)^-}\right)^q = \left(\frac{4}{1-\beta^2 A'_1}\right)^q + O(1)$$

et que

$$(A'_q)^- = A'_q + O(1).$$

Posons

$$\lambda^- = \beta_-\sqrt{q_2^-}, \qquad \lambda = \beta\sqrt{q_2},$$
$$Y(\lambda) = \lambda z + h, \qquad F(\lambda) = E[(1-th^2(Y(\lambda)))^q].$$

Il est clair que $(A'_q)^- - A'_q = F(\lambda^-) - F(\lambda)$ et que la dérivée de $F$ en $\lambda$ est majorée par

$$\int_R K|z|f(z)\,dz \leq K,$$

où $f(z)$ est la densité de la gaussienne standard et $K$ est un majorant de la valeur absolue de la dérivée de la fonction de $Y$ $(1-th^2(Y))^q$.

Il s'ensuit que

$$\forall q \qquad |(A'_q)^- - A'_q| \leq K|\lambda - \lambda^-|.$$

Comme

$$\lambda - \lambda^- = \beta\left(\sqrt{q_2} - \sqrt{\frac{(N-1)q_2^-}{N}}\right) = \beta\frac{(q_2 - q_2^-(N-1)/N)}{\sqrt{q_2} + \sqrt{(N-1)q_2^-/N}},$$

on voit que

$$|\lambda - \lambda^-| \leq K(|q_2 - q_2^-|) = O(2).$$

Il s'ensuit que

$$\forall q \qquad |(A'_q)^- - A'_q| = O(2)$$

et que

$$\beta_-^2 (A'_1)^- = \beta^2 A'_1 + O(2).$$

Il en résulte bien que

$$\left(\frac{4}{1 - \beta_-^2 (A'_1)^-}\right)^q = \left(\frac{4}{1 - \beta^2 A'_1}\right)^q + O(2),$$

et donc le résultat désiré. □

## REFERENCES

[1] AIZENMAN, M., LEBOWITZ, J. L. and RUELLE, D. (1987). Some rigourous results on the Sherrington–Kirkpatrick spin glass model. *Comm. Math. Phys.* **112** 3–20. MR0904135
[2] BARDINA, X., MARQUEZ-CARRERAS, D., ROVIRA, C. and TINDEL, S. (2004). The p-spin interaction model with external field. *Potential Anal.* **21** 311–362. MR2081143
[3] BARDINA, X., MARQUEZ-CARRERAS, D., ROVIRA, C. and TINDEL, S. (2004). Higher order expansions for the overlap of the SK model. In *Seminar on Stochastic Analysis, Random Fields and Applications IV* (R. C. Dalang, M. Dozzi and F. Russo, eds.) 21–43. Birkhäuser, Basel. MR2096278
[4] BOVIER, A., KURKOVA, I. and LÖWE, M. (2002). Fluctuations of the free energy in the REM and the p-spin SK models. *Ann. Probab.* **30** 605–651. MR1905853
[5] COMETS, F. and NEVEU, J. (1995). The Sherrington–Kirkpatrick model of spin glasses and stochastic calculus: The high temperature case. *Comm. Math. Phys.* **166** 549–564. MR1312435
[6] GUERRA, F. and TONINELLI, F. L. (2002). Central limit theorem for fluctuations in the high temperature region of the Sherrington–Kirkpatrick spin glass model. *J. Math. Phys.* **43** 6224–6237. MR1939641

[7] HANEN, A. (2006). Un théorème limite pour les covariances des spins en deux sites dans le modèle de Sherrington–Kirkpatrick avec champ externe. *C. R. Math. Acad. Sci. Paris* **342** 147–150. MR2193663
[8] KURKOVA, I. (2005). Fluctuations of the free energy and overlaps in the high-temperature p-spin SK and Hopfield models. *Markov Process. Related Fields* **11** 55–80. MR2133344
[9] MEZARD, M., PARISI, G. and VIRASORO, M. A. (1987). *Spin Glass Theory and Beyond.* Word Scientific, Teaneck, NJ. MR1026102
[10] TALAGRAND, M. (2003). *Spin Glasses*: *A Challenge for Mathematicians*. Springer, Berlin. MR1993891

69 RUE BARRAULT
PARIS 75013
FRANCE
E-MAIL: AHanen3752@aol.com